Direct sum condition and Artin Approximation in Banach spaces

Matthias Stiefenhofer

ABSTRACT. The system of undetermined coefficients of a bifurcation problem $G[z] = 0$ in Banach spaces is investigated for proving the existence of families of solution curves by use of the implicit function theorem. The main theorem represents an Artin-Tougeron type result in the sense that approximation of order $2k$ ensures exact solutions agreeing up to order $k$ with the approximation [13], [22]. Alternatively, it may be interpreted as Hensel's Lemma in Banach spaces.

In the spirit of [9] and [18], the required surjectivity condition is interpreted as a direct sum condition of order $k$ that allows for solving the remainder equation with respect to graded subspaces derived from an appropriate filtration [24], [25]. In the direction of these subspaces, the determinant can be calculated in a finite dimensional setting, enabling the investigation of secondary global bifurcation phenomena by sign change of Brouwer's degree [18]. The direct sum condition seems to be a generalization of the direct sum condition introduced in [9].

The implicit function theorem delivers stability of $k$ leading coefficients $[\bar{z}_1, \ldots, \bar{z}_k]$ with respect to perturbations of order $2k + 1$ and uniqueness in pointed wedges around the solution curves. Further, denoting by $\bar{X}_\infty$ the subset of arc space $X_\infty$ with $[\bar{z}_1, \ldots, \bar{z}_k]$ fixed, $\bar{X}_\infty$ is iteratively approximated by $\bar{X}_{2k+l}$ as $l \to \infty$. In addition, a lower bound of the Greenberg function $\beta(i)$ of a singularity is constructed by use of a step function obtained from $k$-degree of different solution curves.

Finally, based on Kouchnirenko's theorem [17] and an extension in [6], the results are applied to Newton-polygons where it is shown that the Milnor number of a singularity can be calculated by the sum of $k$-degrees of corresponding solution curves. Simple $ADE$-singularities are investigated in detail.

The main theorem represents a version of strong implicit function theorem in Banach spaces, possibly comparable to theorems in [2]. Moreover, our aim is to extend the direct sum condition of order $k$ from [9] to certain topics in singularity and approximation theory. Generalizations to modules as well as specializations to polynomials seem to be promising.

Contents

*1. Introduction and Notation*

An analytic equation

$$G[z] = 0 \quad \text{with} \quad G[0] = 0, \quad G : B \to \bar{B}$$

$B, \bar{B}$ real or complex Banach spaces is considered, with the aim of finding solution curves $G[z(\varepsilon)] = 0$ through the origin with $\varepsilon \in \mathbb{K} = \mathbb{R}, \mathbb{C}$, $|\varepsilon| \ll 1$. Let us denote by $X_n, n \geq 1$ the set of solutions of the necessary conditions

$$\begin{aligned} T^1[\,z_1\,] &= 0 \\ &\vdots \\ T^n[\,z_1, \ldots, z_n\,] &= 0 \end{aligned} \quad (1)$$

obtained by $\varepsilon$-expansion of

$$G\left[\varepsilon \cdot z_1 + \cdots + \frac{1}{n!}\varepsilon^n \cdot z_n\right] = \sum_{i=1}^n \frac{1}{i!}\varepsilon^i \cdot T^i[\,z_1, \ldots, z_i\,] + \varepsilon^{n+1} \cdot r[\,\varepsilon, z_1, \ldots, z_n\,] = 0$$

with analytic remainder map $r[\cdot]$ and $X_n$ representing the coefficients of approximate solution curves satisfying $G[z(\varepsilon)] = 0$ up to order $O(|\varepsilon|^{n+1})$. If a sequence $z_1, z_2, \ldots$ of coefficients satisfies (1) with $n$ up to infinity (formal power series solution), then the sequence is said to lie in $X_\infty$.

Our aim is to increase $n$ step-by-step up to $n = 2k$ with $k \geq 1$ sufficiently high, yielding initial approximations from $X_{2k}$ that can be lifted to exact solution curves of $G[z] = 0$ by use of the implicit function theorem.

As a preliminary, define the Fréchet derivative $S_1 := G'[0] \in L[B, \bar{B}]$ with null space $N_1 := N[S_1] \subset B$, range $R_1 := R[S_1] \subset \bar{B}$ and corresponding direct sum decompositions of $B$ and $\bar{B}$ according to

$$
\begin{array}{c}
B = N_1^c \oplus N_1 \\
\uparrow \\
\boxed{S_1} \\
\downarrow \\
\bar{B} = R_1 \oplus R_1^c
\end{array}
\qquad (2)
$$

The complementary subspaces $N_1^c$ and $R_1^c$ as well as $R_1$ by itself are assumed to be closed subspaces with continuous bijection $S_1$ between $N_1^c$ and $R_1$ as indicated in (2) by arrows.

Now in case of $S_1$ to be surjective, we obtain $R_1^c = \{0\}$ and the implicit function theorem ensures all solutions of $G[z] = 0$ in the vicinity of $z = 0$ parametrized by $N_1$. In case of $S_1$ not to be surjective, i.e. $R_1^c \neq \{0\}$, we enter the iteration shown in figure 1.

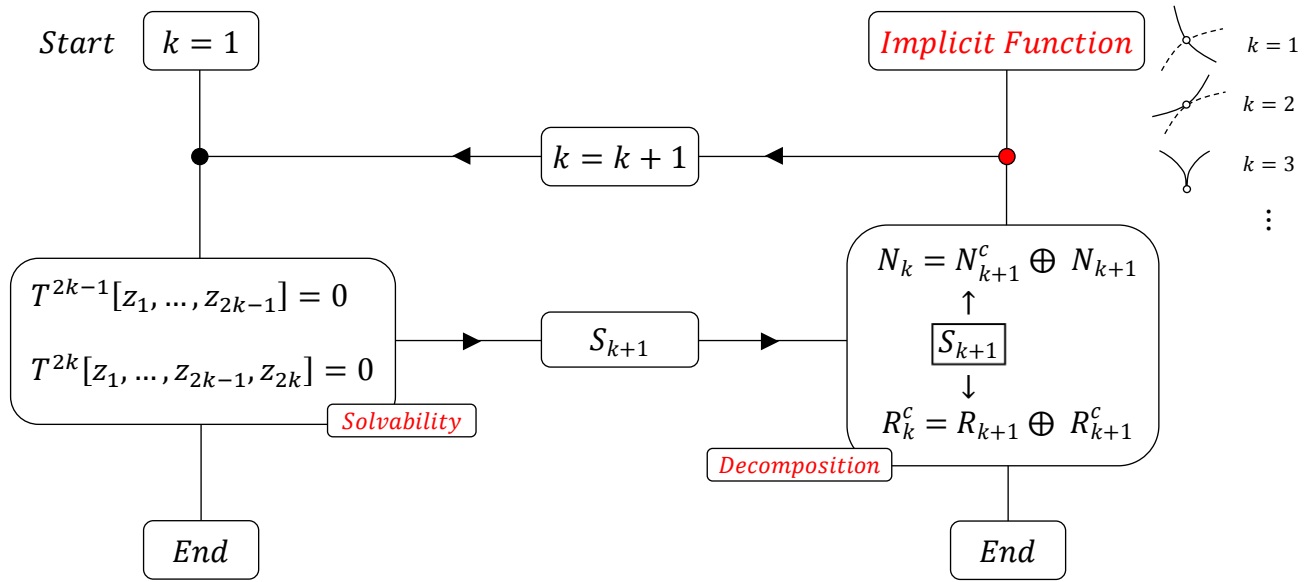

Figure 1 : The iteration process.

The iteration starts top-left with $k = 1$, looking first at the two equations $T^1[z_1] = 0$ and $T^2[z_1, z_2] = 0$. We choose $\bar{z}_1 \in B$ appropriately and restrict $z_2 \in B$ in such a way that for fixed $\bar{z}_1$ all remaining solutions of the two equations

$$T^1[\bar{z}_1] = 0 \quad and \quad T^2[\bar{z}_1, z_2] = 0 \qquad (3)$$

are obtained. Then the linear operator

$$S_2 := P_{R_1^c} \circ 2G''[0]\bar{z}_1 |_{N_1} \in L[N_1, R_1^c] \qquad (4)$$

is defined with continuous projection $P_{R_1^c}$ evaluated with respect to decomposition (2). Note that $S_2$ is mapping from the kernel $N_1$ of $S_1$ into the complement $R_1^c$ of the image $R_1$ of $S_1$, i.e. the second operator $S_2$ is just creating values in the subspace $R_1^c$ that is not reached by the first map-ping. As for $S_1$, we assume for $S_2$ a decomposition of $N_1$ and $R_1^c$ by closed subspaces using null space $N_2 := N[S_2] \subset N_1$ and range $R_2 := R[S_2] \subset R_1^c$ (cf. center-right box in figure 1 with $k = 1$).

Now, if $R_1^c$ already agrees with $R_2$, i.e. $R_2^c = \{0\}$, then the associated sum operator

$$[\,S_1\ \ S_2\,]\begin{pmatrix} n_0 \\ n_1 \end{pmatrix} = S_1 n_0 + S_2 n_1 \in L[\,N_0 \times N_1, \bar{B}\,], \quad N_0 := B \tag{5}$$

is surjective and the implicit function theorem can be applied to an appropriate remainder equation, finally yielding a solution curve of the form

$$z(\varepsilon) = \varepsilon \cdot \bar{z}_1 + \tfrac{1}{2}\varepsilon^2 \cdot z_2(\varepsilon)\,.$$

By choosing different values of $\bar{z}_1$ in (3), other solution curves may be found, as indicated in figure 1 top-right with $k = 1$. The result agrees essentially with the classical bifurcation theorem of simple bifurcation points [7], [8].

If $R_1^c$ is not filled up completely by $R_2$, i.e. $R_2^c \neq \{0\}$, then the iteration parameter $k$ is increased from $k = 1$ to $k = 2$ and the second round in figure 1 gets started with solvability conditions $T^3[\bar{z}_1, z_2, z_3] = 0$ and $T^4[\bar{z}_1, z_2, z_3, z_4] = 0$. These two conditions are used to fix the second coefficient $\bar{z}_2$ appropriately and to restrict $z_3$, $z_4$ to affine subspaces of $B$ in such a way that all remaining solutions of

$$T^3[\,\bar{z}_1, \bar{z}_2, z_3\,] = 0 \qquad and \qquad T^4[\,\bar{z}_1, \bar{z}_2, z_3, z_4\,] = 0$$

are obtained. If no solution at all of $T^3[\bar{z}_1, z_2, z_3] = 0$ and $T^4[\bar{z}_1, z_2, z_3, z_4] = 0$ can be found, then the iteration stops and we have to go back to $k = 1$, possibly restarting the process by use of another value of $\bar{z}_1$.

If an appropriate $\bar{z}_2$ is found, then, analogously to the previous round, the linear operator $S_3 \in L[N_2, R_2^c]$ is defined just creating values in the subspace $R_2^c$ that is not reached by the first two mappings $S_1$ and $S_2$. Further, if decomposition holds with respect to $S_3$ and $N_2$, $R_2^c$, and if the new sum operator

$$[\,S_1\ \ S_2\ \ S_3\,]\begin{pmatrix} n_0 \\ n_1 \\ n_2 \end{pmatrix} = S_1 n_0 + S_2 n_1 + S_3 n_2 \in L[\,N_0 \times N_1 \times N_2, \bar{B}\,] \tag{6}$$

is surjective, we obtain solution curves of the form

$$z(\varepsilon) = \varepsilon \cdot \bar{z}_1 + \tfrac{1}{2}\varepsilon^2 \cdot \bar{z}_2 + \tfrac{1}{6}\varepsilon^3 \cdot z_3(\varepsilon) \tag{7}$$

dealing with bifurcation points of tangentially touching branches with different curvature, as indicated in figure 1 top-right with $k = 2$. The result can be found in [15], [16].

Now in this paper, we show in a constructive way, how to continue the iteration process arbitrarily, thereby enlarging the range of the sum operator step-by-step, hopefully until surjectivity of the sum operator is reached and application of the implicit function theorem is possible.

For general $k \geq 1$ the solvability conditions $T^{2k-1} = 0$ and $T^{2k} = 0$ are used to fix the coefficient $\bar{z}_k$ and to restrict the remaining free coefficients $z_{k+1}, \ldots, z_{2k}$ in such a way that all solutions of

$$T^{2k-1}[\,\bar{z}_1, \ldots, \bar{z}_k, z_{k+1}, \ldots, z_{2k-1}\,] = 0 \quad and \quad T^{2k}[\,\bar{z}_1, \ldots, \bar{z}_k, z_{k+1}, \ldots, z_{2k}\,] = 0 \tag{8}$$

are found. In some more detail, with fixed $k$-tuple $[\bar{z}_1, \ldots, \bar{z}_k]$, the $k$ free variables $z_{k+1}, \ldots, z_{2k}$ can exactly be chosen from an affine subspace in $B^k$ for satisfying (8), thus defining the following subset of $X_{2k}$

$$\bar{X}_{2k} := \{\,[\,z_1, \ldots, z_k, z_{k+1}, \ldots, z_{2k}\,] \in X_{2k}\ |\ z_1 = \bar{z}_1, \ldots, z_k = \bar{z}_k\,\}\,.$$

In general, $\bar{X}_n, n \geq k$ defines the subset in $X_n$ with fixed leading coefficients $[\bar{z}_1, \ldots, \bar{z}_k] \in B^k$. The basis for the effective construction of $\bar{X}_{2k}$ is given by some sort of filtration in $B$

$$B = N_0 \supset N_1 \supset \cdots \supset N_k \supset N_{k+1} \tag{9}$$

that is sequentially defined by the kernels of the linear operators $S_1, \ldots, S_{k+1}$. In addition, the filtration can be used to obtain direct sum decompositions of $B$ and $\bar{B}$ by graded subspaces according to

$$\begin{array}{c} B = N_1^c \oplus N_2^c \oplus \cdots \oplus N_{k+1}^c \oplus N_{k+1} \\ \uparrow \quad \uparrow \quad \quad \uparrow \\ \boxed{S_1} \quad \boxed{S_2} \quad \quad \boxed{S_{k+1}} \\ \downarrow \quad \downarrow \quad \quad \downarrow \\ \bar{B} = R_1 \oplus R_2 \oplus \cdots \oplus R_{k+1} \oplus R_{k+1}^c \end{array} \tag{10}$$

Then, if the associated sum operator

$$[\, S_1 \cdots S_{k+1}\, ] \begin{pmatrix} n_0 \\ \vdots \\ n_k \end{pmatrix} = S_1 n_0 + \cdots + S_{k+1} n_k \in L[\, N_0 \times \cdots \times N_k, \bar{B}\, ] \tag{11}$$

is surjective, i.e. $R_{k+1}^c = \{0\}$ in (10), we obtain solution curves of the form

$$z(\varepsilon) = \varepsilon \bar{z}_1 + \cdots + \frac{1}{k!}\varepsilon^k \bar{z}_k + \frac{1}{(k+1)!}\varepsilon^{k+1} z_{k+1}(\varepsilon)\,.$$

Following [15], we call a $k$-tuple $[\bar{z}_1, \ldots, \bar{z}_k]$ regular of degree $k$, if it represents the leading coefficients of an element in $\bar{X}_{2k}$ with corresponding sum operator (11) showing surjectivity. Finally, define for $X_n$ (and $\bar{X}_n$) the truncation maps $\pi_l: X_n \to X_l, n \geq l$ by deleting the components $z_{l+1}, \ldots, z_n$.

Before stating the results in detail within section 3, we aim to present several examples in section 2 for showing the working principle of the iteration with respect to different types of applications. The example section, which is somewhat long, may be skipped. Applications to singularly perturbed systems in case of $k = 1$ may be found in [26], [27] and [28]. A blow-up prodedure with respect to homoclinic orbits is performed in [3].

## 2. Examples

Example 1 outlines in detail the building process of the sum operator by reference to a simple example with $G: \mathbb{R}^3 \to \mathbb{R}^2$. In addition, the stability and uniqueness results from section 3 are addressed and a solution curve that cannot be established with the iteration process is presented, based on the existence of nonisolated critical points within the singular locus $X_{sing}$.

Example 2 is dealing with $ADE$-singularities [1], [4] and $k$-degree of their solution curves. Additionally, some formulas concerning the operators $S_{k+1}$ are given in the simplifying cases of $G: \mathbb{K}^2 \to \mathbb{K}$ as well as an existing trivial solution curve $y \equiv 0$. Finally, the Milnor number of $ADE$-singularities is calculated by $k$-degree of their solution curves. This result is a special application of Corollary 4 of section 4 which allows to calculate the Milnor number by some characteristics of the Newton-polygon and $k$-degree of corresponding solution curves. The proof of Corollary 4 is based on Kouchnirenko's planar theorem [1], [17].

In Example 3, we present some basic relations between the direct sum condition from local/global bifurcation theory [18], [19] and decomposition (10). Here, the behavior of the determinant with respect to a complementary subspace of the solution curve is calculated and used to consider possible sign change of Brouwer's degree assuring bifurcation of secondary solutions. Decomposition (10) seems to be a generalization of the direct sum condition introduced by J. López-Gómez in [9].

Finally, examples 4 and 5 are dealing with families of solution curves, creating 2- and 3-dimensional surfaces of solutions. Example 4 looks at the Whitney umbrella under perturbations and examines $k$-degree and stability of different segments of the umbrella.

Example 5 revisits an example of G. Belitskii and D. Kerner. In [2], the example illustrates the advantages of the strong implicit function theorem for modules developed in [2], compared to the implicit function theorems of Tougeron [30] and Fisher [11]. By application of the iteration from figure 1 to the example, we ascertain $k$-degree of the solutions, implying stability with respect to $(x, y)$-perturbations of order $2k + 1$. The result may be an extension to [2], since the perturbations are also allowed to depend on $y$.

**Example 1 – Sum Operator :** Consider the equation

$$G[\,z\,] = G[\,x_1, x_2, x_3\,] = \begin{pmatrix} x_1 x_3 \\ x_2^2 - x_2 x_1^2 \end{pmatrix} = \begin{pmatrix} 0 \\ 0 \end{pmatrix} \quad \text{with} \quad G : \mathbb{R}^3 \to \mathbb{R}^2 \qquad (12)$$

and solutions obviously given by $x_1$-axis, $x_3$-axis and the parabola $x_2 = x_1^2$ positioned within $(x_1, x_2)$-plane, as shown in figure 2 top left.

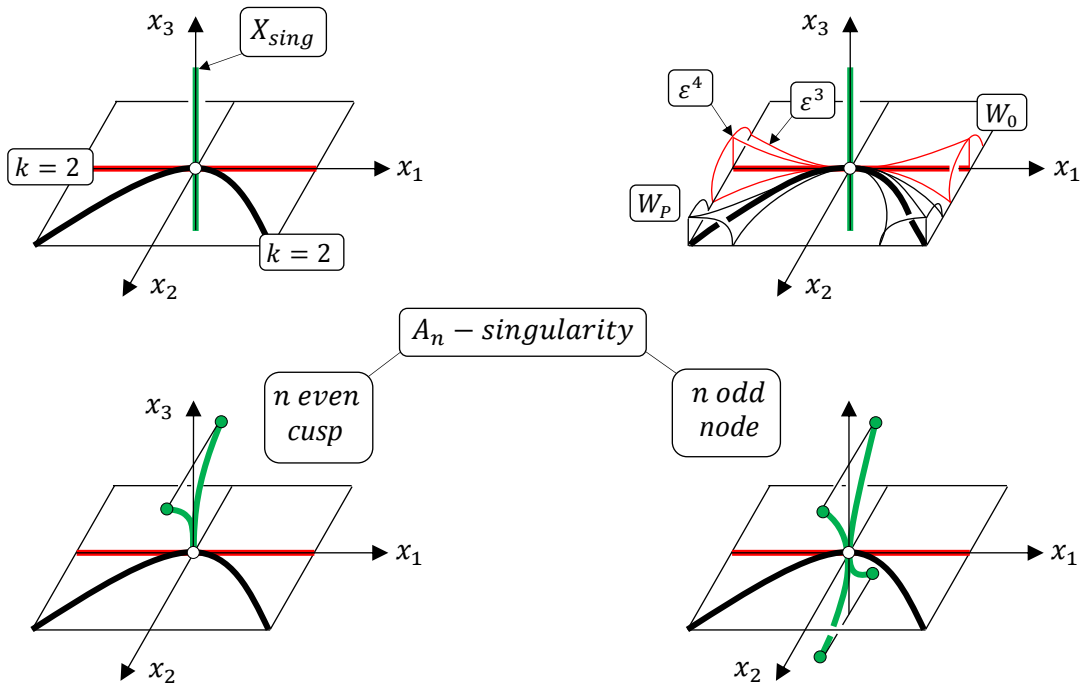

Figure 2 : Basic solution curves (top left) with $k$-degree of regularity.

Pointed wedge-like neighborhoods $W_0$ and $W_P$ of uniqueness (top right).

$A_n$-singularities under perturbation (cusp left, node right).

Let us now perform the iteration procedure step by step, aiming to prove the existence of the different solution curves by building up the corresponding sum operators. First and according to (2), define $S_1 := G'[0] = 0 \in L[\mathbb{R}^3, \mathbb{R}^2]$ with null space $N_1 = \mathbb{R}^3$, range $R_1 = \{0\}$ and belonging complements $N_1^c = \{0\}$ as well as $R_1^c = \mathbb{R}^2$, i.e. the first linear operator $S_1$ is completely singular, hence not contributing to the surjectivity of the sum operator.

Next, using the abbreviations $z_i = (u_i, v_i, w_i)^T \in \mathbb{R}^3$ and $G_0^i = \frac{d^i}{dz^i} G[0] \in L^i[B, \bar{B}]$ for $= 1, 2, \ldots$, the iteration is started with $k = 1$ by the first two equations

$$T^1[\, z_1 \,] = \overbrace{G_0^1}^{=0} z_1 = 0 \tag{13}$$

$$T^2[\, z_1, z_2 \,] = \underbrace{G_0^1}_{=0} z_2 + G_0^2 z_1^2 = \begin{pmatrix} 2u_1 w_1 \\ 2v_1^2 \end{pmatrix} = 0$$

with solutions $z_1 = (u_1, 0, 0)^T, u_1 \in \mathbb{R}$ and $z_1 = (0, 0, w_1)^T, w_1 \in \mathbb{R}$, as well as $z_2 = (u_2, v_2, w_2)^T \in \mathbb{R}^3$ in both cases. First, we concentrate on the red marked solution curve in figure 1 along the $x_1$-axis by setting $\bar{z}_1 = (1, 0, 0)^T$. The case $\bar{z}_1 = (0, 0, 1)^T$ is treated at the end of the example. Then, under consideration of $R_1^c = \mathbb{R}^2$, $N_1 = \mathbb{R}^3$ and (4), the second linear operator reads

$$S_2 := P_{R_1^c} \circ 2G_0^2 \bar{z}_1 \,|_{N_1} = 2G_0^2 \bar{z}_1 = 2\begin{pmatrix} 0 & 0 & 1 \\ 0 & 0 & 0 \end{pmatrix}$$

with $N_2 = \{(1,0,0)^T, (0,1,0)^T\}$ and $R_2 = \{(1,0)^T\}$, yielding at this stage of the process the decompositions

$$B = \mathbb{R}^3 = \overbrace{\left\{\begin{pmatrix}0\\0\\0\end{pmatrix}\right\}}^{=N_1^c} \oplus \overbrace{\left\{\begin{pmatrix}0\\0\\1\end{pmatrix}\right\}}^{=N_2^c} \oplus \overbrace{\left\{\begin{pmatrix}1\\0\\0\end{pmatrix}, \begin{pmatrix}0\\1\\0\end{pmatrix}\right\}}^{=N_2}$$

$$\uparrow \qquad \uparrow$$
$$\boxed{S_1 = 0} \quad \boxed{S_2 \neq 0} \tag{14}$$
$$\downarrow \qquad \downarrow$$

$$\bar{B} = \mathbb{R}^2 = \underbrace{\left\{\begin{pmatrix}0\\0\end{pmatrix}\right\}}_{=R_1} \oplus \underbrace{\left\{\begin{pmatrix}1\\0\end{pmatrix}\right\}}_{=R_2} \oplus \underbrace{\left\{\begin{pmatrix}0\\1\end{pmatrix}\right\}}_{=R_2^c}$$

Here, brackets $\{\cdots\}$ denote the subspace spanned by the elements within the brackets. Obviously, the sum operator from (5) is not yet surjective and $k$ has to be increased from $k = 1$ to $k = 2$ according to figure 1. Then, the solvability condition reads

$$T^3[\, \bar{z}_1, z_2, z_3 \,] = \overbrace{G_0^1}^{=0} z_3 + 3G_0^2 \bar{z}_1 z_2 + \overbrace{G_0^3 \bar{z}_1^3}^{=0} = 0$$

$$T^4[\, \bar{z}_1, z_2, z_3, z_4 \,] = \underbrace{G_0^1}_{=0} z_4 + 4G_0^2 \bar{z}_1 z_3 + 3G_0^2 z_2^2 + 6G_0^3 \bar{z}_1^2 z_2 + \underbrace{G_0^4\, \bar{z}_1^4}_{=0} = 0$$

or equivalently

$$3\begin{pmatrix}w_2\\0\end{pmatrix} = \begin{pmatrix}0\\0\end{pmatrix} \quad \wedge \quad 4\begin{pmatrix}w_3\\0\end{pmatrix} + 3\begin{pmatrix}0\\2v_2^2\end{pmatrix} + 6\begin{pmatrix}0\\-2v_2\end{pmatrix} = \begin{pmatrix}0\\0\end{pmatrix}$$

yielding the two sets of solutions for $z_2, z_3$ and $z_4$

$$z_2 = (u_2, 0, 0)^T, u_2 \in \mathbb{R}, \quad z_3 = (u_3, v_3, 0)^T, (u_3, v_3) \in \mathbb{R}^2, \quad z_4 \in \mathbb{R}^3$$
$$z_2 = (u_2, 2, 0)^T, u_2 \in \mathbb{R}, \quad z_3 = (u_3, v_3, 0)^T, (u_3, v_3) \in \mathbb{R}^2, \quad z_4 \in \mathbb{R}^3 \tag{15}$$

First, we concentrate on the upper case, fixing $z_2$ by $\bar{z}_2 = (1, 0, 0)^T$ for definition of the linear map

$$S_3 := P_{R_2^c} \circ \left[\, 6G_0^2 \bar{z}_2 + 6G_0^3 \bar{z}_1^2 - 12G_0^2 \bar{z}_1 \circ S_1^{-1} P_{R_1} \circ G_0^2 \bar{z}_1 \,\right]_{|N_2} \in L[\, N_2, R_2^c \,] \tag{16}$$

with continuous projections $P_{R_2^c}$ and $P_{R_1}$ evaluated with respect to decomposition (14). Then, by (16) we obtain

$$S_3 = 6 \cdot P_{R_2^c} \circ \begin{pmatrix} 0 & 0 & 1 \\ 0 & -2 & 0 \end{pmatrix}_{|N_2} = 6 \cdot \begin{pmatrix} 0 & 0 & 0 \\ 0 & -2 & 0 \end{pmatrix}$$

implying $N_3 = \{(1,0,0)^T\} \subset N_2$ and $R_3 = \{(0,1)^T\} = R_2^c$ as well as the final decomposition

$$B = \mathbb{R}^3 = \overbrace{\left\{\begin{pmatrix} 0 \\ 0 \\ 0 \end{pmatrix}\right\}}^{=N_1^c} \oplus \overbrace{\left\{\begin{pmatrix} 0 \\ 0 \\ 1 \end{pmatrix}\right\}}^{=N_2^c} \oplus \overbrace{\left\{\begin{pmatrix} 0 \\ 1 \\ 0 \end{pmatrix}\right\}}^{=N_3^c} \oplus \overbrace{\left\{\begin{pmatrix} 1 \\ 0 \\ 0 \end{pmatrix}\right\}}^{=N_3}$$

$$\uparrow \qquad \uparrow \qquad \uparrow$$
$$\boxed{S_1 = 0} \quad \boxed{S_2 \neq 0} \quad \boxed{S_3 \neq 0} \qquad (17)$$
$$\downarrow \qquad \downarrow \qquad \downarrow$$

$$\bar{B} = \mathbb{R}^2 = \underbrace{\left\{\begin{pmatrix} 0 \\ 0 \end{pmatrix}\right\}}_{=R_1} \oplus \underbrace{\left\{\begin{pmatrix} 1 \\ 0 \end{pmatrix}\right\}}_{=R_2} \oplus \underbrace{\left\{\begin{pmatrix} 0 \\ 1 \end{pmatrix}\right\}}_{=R_3\, =R_2^c}$$

with surjectivity of the sum operator $[S_1 \ S_2 \ S_3]$ from (6). Hence, according to the implicit function theorem and (7), a solution curve of the form

$$z_0(\varepsilon) = \varepsilon \cdot \underbrace{\begin{pmatrix} 1 \\ 0 \\ 0 \end{pmatrix}}_{=\bar{z}_1} + \tfrac{1}{2}\varepsilon^2 \cdot \underbrace{\begin{pmatrix} 1 \\ 0 \\ 0 \end{pmatrix}}_{=\bar{z}_2} + O(|\varepsilon|^3) \qquad (18)$$

exists, showing that the 2-tuple $[\bar{z}_1, \bar{z}_2]$ is regular of degree $k = 2$. This solution curve represents the red marked solutions along the $x_1$-axis in figure 2.

In the next step, the black parabola in figure 2 is established by choosing the lower case in (15) with $u_2 = 1$, $\bar{z}_2 = (1,2,0)^T$, yielding again decomposition (17) and parametrization of the parabola by

$$z_P(\varepsilon) = \varepsilon \cdot \underbrace{\begin{pmatrix} 1 \\ 0 \\ 0 \end{pmatrix}}_{=\bar{z}_1} + \tfrac{1}{2}\varepsilon^2 \cdot \underbrace{\begin{pmatrix} 1 \\ 2 \\ 0 \end{pmatrix}}_{=\bar{z}_2} + O(|\varepsilon|^3). \qquad (19)$$

Note also that choosing other values in (15) with respect to $u_2$, e.g. $u_2 = 0$, delivers same solution orbits, alternatively parametrized by $\bar{z}_2 = (0,0,0)^T$ in (18) and $\bar{z}_2 = (0,2,0)^T$ in (19).

Both of the solution curves (18) and (19) are regular of degree $k = 2$, where in the construction of the leading coefficients $\bar{z}_1$ and $\bar{z}_2$ the derivatives $G_0^1, \ldots, G_0^{2k} = G_0^4$ are involved. The sum operator $[S_1 \ S_2 \ S_3]$ is affected by $\bar{z}_1, \bar{z}_k = \bar{z}_2$ and derivatives $G_0^1, \ldots, G_0^{k+1} = G_0^3$. These principal dependencies carry over to general $k \geq 1$.

In Corollary 1 of section 3, we will see that the following stability results are valid with respect to solution curves (18) and (19). Perturbations of $G[z]$ of order $O(|z|^{2k+1}) = O(|z|^5)$ do not change the two leading coefficients $\bar{z}_1$ and $\bar{z}_2$, only varying higher order terms $O(|\varepsilon|^3)$. Secondly, perturbations of $G[z]$ of the form $G[z] + c \cdot H[z]$ with $H[z] = O(|z|^{2k}) = O(|z|^4)$ will not destroy the solution curves if the constant $c$ is chosen sufficiently small ($\bar{z}_1$ remains unchanged, $\bar{z}_2$ may vary). Finally, perturbations of order $O(|z|^{2k-1}) = O(|z|^3)$ may destroy the solution curves.

For example, perturbing $G[z]$ by $H[z] = O(|z|^4)$ according to

$$G[\,z\,] + c \cdot H[\,z\,] = \begin{pmatrix} x_1 x_3 \\ x_2^2 - x_2 x_1^2 \end{pmatrix} + c \cdot \begin{pmatrix} 0 \\ x_1^4 \end{pmatrix}$$

will not destroy the two solution curves for $c < 0.25$, whereas at $c = 0.25$, the two solution curves merge and disappear for $c > 0.25$.

In Corollary 3 of section 3, it is additionally seen that the implicit function theorem delivers uniqueness in the following sense. Solution curves with finite $k$-degree of regularity are embedded in a pointed wedge-like neighborhood within $B$ that shrinks to zero with different orders of $\varepsilon$ with respect to the decomposition of $B$ in (10), (17). Concerning the solution curves (18) and (19) with $k = 2$, the neighborhoods $W_0$ of $z_0(\varepsilon)$ and $W_P$ of $z_P(\varepsilon)$ are given by

$$W_{0/P}: \; z_{0/P}(\varepsilon) + \varepsilon^3 \cdot \overbrace{\begin{pmatrix} 0 \\ x_2 \\ 0 \end{pmatrix}}^{\subset N_3^c} + \varepsilon^4 \cdot \overbrace{\begin{pmatrix} 0 \\ 0 \\ x_3 \end{pmatrix}}^{\subset N_2^c} \quad with \quad |x_2|, |x_3| \ll 1$$

as indicated in figure 2 top-right. Hence, in the directions of the subspaces $N_3^c$ ($x_2$-axis) and $N_2^c$ ($x_3$-axis), the wedges are shrinking to zero by order of $O(|\varepsilon|^3)$ and $O(|\varepsilon|^4)$ respectively.

It remains to look at the solutions along the $x_3$-axis by choice of $\bar{z}_1 = (0,0,1)^T$ in (13), implying a decomposition after the second round according to

$$B = \mathbb{R}^3 = \overbrace{\left\{\begin{pmatrix} 0 \\ 0 \\ 0 \end{pmatrix}\right\}}^{=N_1^c} \oplus \overbrace{\left\{\begin{pmatrix} 1 \\ 0 \\ 0 \end{pmatrix}\right\}}^{=N_2^c} \oplus \overbrace{\left\{\begin{pmatrix} 0 \\ 0 \\ 0 \end{pmatrix}\right\}}^{=N_3^c} \oplus \overbrace{\left\{\begin{pmatrix} 0 \\ 1 \\ 0 \end{pmatrix}, \begin{pmatrix} 0 \\ 0 \\ 1 \end{pmatrix}\right\}}^{=N_3}$$

$$\uparrow \qquad \uparrow \qquad \uparrow$$
$$\boxed{S_1 = 0} \quad \boxed{S_2 \neq 0} \quad \boxed{S_3 = 0}$$
$$\downarrow \qquad \downarrow \qquad \downarrow$$

$$\bar{B} = \mathbb{R}^2 = \underbrace{\left\{\begin{pmatrix} 0 \\ 0 \end{pmatrix}\right\}}_{=R_1} \oplus \underbrace{\left\{\begin{pmatrix} 1 \\ 0 \end{pmatrix}\right\}}_{=R_2} \oplus \underbrace{\left\{\begin{pmatrix} 0 \\ 0 \end{pmatrix}\right\}}_{=R_3} \oplus \underbrace{\left\{\begin{pmatrix} 0 \\ 1 \end{pmatrix}\right\}}_{=R_3^c}$$

Obviously by $R_3^c \neq \{0\}$, the sum operator $[S_1 \; S_2 \; S_3]$ is not surjective. Now, it is easy to see from (12) that the iteration yields $S_i = 0$ for $i \geq 3$, and hence it is not possible to increase the range of the sum operator any more for proving the existence of the solutions along the $x_3$-axis. This is not astonishing, since the solution curve along the $x_3$-axis is not stable with respect to perturbations of $G[z]$ of arbitrary high order. For example, perturbing $G[z]$ by

$$G[\,x_1, x_2, x_3\,] - \begin{pmatrix} 0 \\ x_3^{n+1} \end{pmatrix} = \begin{pmatrix} x_1 x_3 \\ x_2^2 - x_2 x_1^2 - x_3^{n+1} \end{pmatrix}$$

may destroy or split the solutions along the $x_3$-axis for arbitrary $n \geq 4$, as shown in figure 2 bottom left and right. For $n$ even ($n$ odd), the $x_3$-axis turns into a cusp (node) of the form $x_2^2 - x_3^{n+1} = 0$ lying within $(x_2, x_3)$-plane, i.e. an $A_n$-singularity appears [1], [4]. In case of the cusp, the solutions along the negative $x_3$-axis are completely destroyed.

Note also that the $x_3$-axis represents the singular locus $X_{sing}$ of (12), consisting of nonisolated critical points with $G'[0,0,x_3] = 0$.

In the next example, we are considering $A_n$-singularities in $\mathbb{R}^2$ in some more detail, i.e. simple $ADE$-singularities are analyzed with respect to $k$-degree of their solution curves. Moreover, a close relation between the Milnor number $\mu$ and $k$-degree is established.

**Example 2 – Milnor Number and *ADE*-Singularities :** First note that in case of $G: \mathbb{K}^2 \to \mathbb{K}$, the sum operator reaches surjectivity as soon as the first operator $S_{k+1} \in L[\mathbb{K}^2, \mathbb{K}]$ occurs that is different from the zero operator, i.e. surjectivity is reached in one step by

$$B = \mathbb{K}^2 = \overbrace{\left\{\begin{pmatrix}0\\0\end{pmatrix}\right\}}^{=N_1^c} \oplus \cdots \oplus \overbrace{\left\{\begin{pmatrix}0\\0\end{pmatrix}\right\}}^{=N_k^c} \oplus N_{k+1}^c \oplus N_{k+1}$$

$$\uparrow \qquad\qquad \uparrow \qquad \uparrow$$
$$\boxed{S_1 = 0} \quad \boxed{S_k = 0} \; \boxed{\color{red}S_{k+1} \neq 0}$$
$$\downarrow \qquad\qquad \downarrow \qquad \downarrow$$

$$\bar{B} = \mathbb{K} \;=\; \underbrace{\{0\}}_{=R_1} \oplus \cdots \oplus \underbrace{\{0\}}_{=R_k} \oplus \underbrace{\{1\}}_{=R_{k+1}}$$

Exemplarily, we start the investigation of *ADE*-singularities with a real $D_5$-singularity given by

$$G[z] = G[x, y] = x^2 y - y^4 \quad \text{with} \quad G : \mathbb{R}^2 \to \mathbb{R}$$

and corresponding basic solution curve along the $x$-axis

$$z_1(\varepsilon) = \varepsilon \cdot \underbrace{\begin{pmatrix}1\\0\end{pmatrix}}_{=\bar{z}_1} \quad \text{with} \quad \bar{z}_i = \begin{pmatrix}0\\0\end{pmatrix}, \; i \geq 2 \tag{20}$$

as well as a cusp curve emanating in the direction of the $y$-axis according to

$$z_2(\varepsilon) = \tfrac{1}{2}\varepsilon^2 \cdot \underbrace{\begin{pmatrix}0\\2\end{pmatrix}}_{=\bar{z}_2} + \tfrac{1}{6}\varepsilon^3 \cdot \underbrace{\begin{pmatrix}6\\0\end{pmatrix}}_{=\bar{z}_3} \quad \text{with} \quad \bar{z}_1 = \begin{pmatrix}0\\0\end{pmatrix} \quad \text{and} \quad \bar{z}_i = \begin{pmatrix}0\\0\end{pmatrix}, \; i \geq 4. \tag{21}$$

Our aim is to determine $k$-degree of regularity of the two solution curves, where the iteration in figure 1 reduces to check for the first operator $S_{k+1}$ to become different from zero.

Now, for general $G: \mathbb{K}^2 \to \mathbb{K}$ and for $k = 0, \ldots, 5$, the operators $S_{k+1}$ are explicitly given by

(22)

$$S_1 = G_0^1$$

$$S_2 = 2 G_0^2 \bar{z}_1$$

$$S_3 = 6 G_0^2 \bar{z}_2 + 6 G_0^3 \bar{z}_1^2$$

$$S_4 = 20 G_0^2 \bar{z}_3 + 60 G_0^3 \bar{z}_1 \bar{z}_2 + 20 G_0^4 \bar{z}_1^3$$

$$S_5 = 70 G_0^2 \bar{z}_4 + 210 G_0^3 \bar{z}_2^2 + 280 G_0^3 \bar{z}_1 \bar{z}_3 + 420 G_0^4 \bar{z}_1^2 \bar{z}_2 + 70 G_0^5 \bar{z}_1^4$$

$$S_6 = 252 G_0^2 \bar{z}_5 + 2520 G_0^3 \bar{z}_2 \bar{z}_3 + 1260 G_0^3 \bar{z}_1 \bar{z}_4 + 3780 G_0^4 \bar{z}_1 \bar{z}_2^2 + 2520 G_0^4 \bar{z}_1^2 \bar{z}_3 + 2520 G_0^5 \bar{z}_1^3 \bar{z}_2 + 252 G_0^6 \bar{z}_1^5$$

based on the formula

(23)

$$S_{k+1} = \frac{1}{k!} \sum_{\beta=1}^{2k} G_0^\beta \sum_{n_1+\cdots+n_{k-1}+1=\beta \;\wedge\; 1\cdot n_1+\cdots+(k-1)\cdot n_{k-1}=k} \frac{(2k)!}{n_1!\cdots n_{k-1}!} \prod_{\tau=1}^{k-1} \left(\frac{1}{\tau!}\bar{z}_\tau\right)^{n_\tau}$$

Due to $\bar{z}_i = 0, i \geq 2$ and $G_{x^i}^0 = 0, i \geq 1$, formula (23) simplifies heavily in case of an existing trivial solution curve along the $x$-axis according to

$$S_{k+1} = \underbrace{c}_{\neq 0} \cdot \begin{bmatrix} 0 & G_{yx^k}^0 \end{bmatrix} \in L[\,\mathbb{K}^2, \mathbb{K}\,]. \tag{24}$$

The simplicity of (22), (23) and (24), based on $G: \mathbb{K}^2 \to \mathbb{K}$, should not be misleading. The number of summands composing the operator $S_{k+1}$ grows very fast with respect to $k$ when dealing with general Banach spaces $B$ and $\bar{B}$, e.g. in general $S_5$ comprises 97 summands. On the other side, the operators are recursively defined and can easily be obtained by computer. Serious simplification occurs in case of $G: \mathbb{K}^2 \to \mathbb{K}$ according to (22) and (23), in case of an existing trivial solution curve and in case of polynomials.

Now, concerning the sigularity $D_5$ and the trivial solutions $z_1(\varepsilon)$ from (20), the first operator different from zero appears with $k = k_1 = 2$ in the form

$$S_3 = 6 \cdot G_0^3 \bar{z}_1^2 = 6 \cdot \begin{bmatrix} 0 & G_{yx^2}^0 \end{bmatrix} = 6 \cdot [0\ 2] \neq [0\ 0]\ surjective.$$

Concerning the cusp curve (21), the first operator different from zero is given by $k = k_2 = 5$ according to

$$S_6 = 2520 \cdot G_0^3 \bar{z}_2 \bar{z}_3 = 2520 \cdot [24\ 0] \neq [0\ 0]\ surjective.$$

In figure 3, the solution curves of $D_5$ are shown together with their pointed wedges of uniqueness by order of $O(|\varepsilon|^3)$ and $O(|\varepsilon|^6)$ respectively. In addition, the possibility to calculate the Milnor number $\mu = 5$ of $D_5$ by the sum of $k$-degrees is indicated on the left hand side, where $ord(G)$ denotes lowest degree within $G[x, y]$.

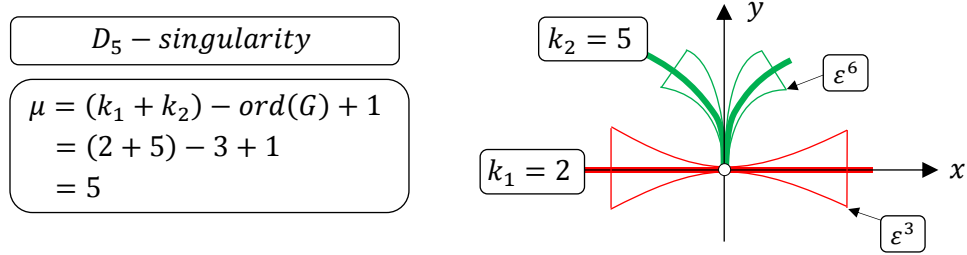

Figure 3 : $D_5$-singularity with solution curves and $k$-degree of regularity.

Note also that the singularity $D_5$ is 4-$R$-determined with respect to right equivalence [12] which means that perturbations of order $O(|x, y|^5)$ do not destroy the solution curves, but may change the direction along which the curves emanate from the origin. Now, the stability result of Corollary 1 in section 3 asserts that perturbations of order $2k_1 + 1 = 5$ do not change the two leading terms $\bar{z}_1 = (1,0)^T$ and $\bar{z}_2 = (0,0)^T$ of the basic solution curve $z_1(\varepsilon)$ in (20). Concerning the cusp curve $z_2(\varepsilon)$, we can only state that perturbations of order $2k_2 + 1 = 11$ do not change the leading coefficients $\bar{z}_1, \cdots, \bar{z}_5$ in (21).

Along the same lines as for $D_5$, other singularities may be investigated, finally implying table 1 that summarizes some results of real $ADE$-singularities with a complete set of solution curves.

| Singularity | $G[x,y]$ | Milnor $\mu$ | $k-$degree | solutions |
|---|---|---|---|---|
| $A_{2n-1}, n \geq 1$ | $x^2 - y^{2n}$ | $2n-1$ | $n, n$ | 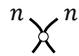 |
| $A_{2n}, n \geq 1$ | $x^2 - y^{2n+1}$ | $2n$ | $2n+1$ | 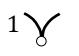 |
| $D_{2n}, n \geq 2$ | $x^2 y - y^{2n-1}$ | $2n$ | $2, n, n$ | 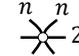 |
| $D_{2n+1}, n \geq 2$ | $x^2 y - y^{2n}$ | $2n+1$ | $2, 2n+1$ | 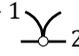 |
| $E_6$ | $x^3 - y^4$ | $6$ | $8$ | 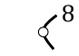 |
| $E_7$ | $x^3 - xy^3$ | $7$ | $3, 6$ | 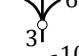 |
| $E_8$ | $x^3 - y^5$ | $8$ | $10$ | 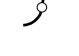 |

Table 1 : Simple $ADE$-singularities with Milnor number $\mu$ and $k$-degree of regularity.

In particular, each of the singularities satisfies

$$\mu = \sum k_i - ord(G) + 1 \,. \tag{25}$$

In case of complex $ADE$-singularities we obtain the same result, where the existence of a complete set of solution curves is assured by algebraic closure.

Within section 4 of this paper, relation (25) is generalized by use of the iteration from figure 1 and Kouchnirenko's planar theorem [17] along the following lines. If $G[x,y]$ defines a convenient Newton-polygon with a complete set of simple zeros for every segment, then the Milnor number $\mu$ can be calculated according to

$$\mu = \sum_{i=1}^{\tau} r_i \cdot k_i - ord(G) + 1 \,. \tag{26}$$

Here $\tau \geq 1$ denotes the number of segments of the Newton-polygon and $r_i \geq 1$ is given by the degree of the homogenous polynomial associated with each segment [6].

**Example 3 - Global Bifurcation and direct sum condition :** In some situations, the direct sum decomposition (10) allows to compute a topological degree along a solution curve. For simplicity, assume the subspace $N^c := N_1^c \oplus N_2^c \oplus \cdots \oplus N_{k+1}^c$ to be finite dimensional and $[\bar{z}_1, \ldots, \bar{z}_k]$ to be $k$-regular with associated decomposition

$$\begin{array}{c} \overbrace{\phantom{N_1^c \oplus N_2^c \oplus \cdots \oplus N_{k+1}^c}}^{=:N^c} \\ B = N_1^c \oplus N_2^c \oplus \cdots \oplus N_{k+1}^c \oplus N_{k+1} \\ \uparrow \quad \uparrow \quad \quad \uparrow \\ \boxed{S_1} \ \boxed{S_2} \quad \boxed{S_{k+1}} \\ \downarrow \quad \downarrow \quad \quad \downarrow \\ \bar{B} = R_1 \oplus R_2 \oplus \cdots \oplus R_{k+1} \end{array} \tag{27}$$

and corresponding solution curve

$$z(\varepsilon) = \tfrac{1}{l!}\varepsilon^l \cdot \underbrace{\bar{z}_l}_{\neq 0} + \cdots + \tfrac{1}{k!}\varepsilon^k \cdot \bar{z}_k + \tfrac{1}{(k+1)!} \cdot \varepsilon^{k+1} z_{k+1}(\varepsilon) \,. \tag{28}$$

Here $\bar{z}_l, l \geq 1$ denotes the first coefficient within $[\bar{z}_1, \ldots, \bar{z}_k]$ different from zero. Then, the linearization of $G[z]$ in the direction of the subspace $N^c$ can be calculated along the solutions $z(\varepsilon)$, where it can be shown that the determinant of $G_{n^c}[z(\varepsilon)] \in L[N^c, \bar{B}]$ is given by an $\varepsilon$-expansion of the form

$$det\{ G_{n^c}[ z(\varepsilon) ] \} = c \cdot \varepsilon^\chi + O(|\varepsilon|^{\chi+1}), \quad c \neq 0 \quad (29)$$

with lowest exponent

$$\chi = 1 \cdot dim\, N_2^c + \cdots + k \cdot dim\, N_{k+1}^c \quad (30)$$

simply composed by the dimensions of $N_2^c, \ldots, N_{k+1}^c$ (or equivalently by the dimensions of $R_2, \ldots, R_{k+1}$). In case of $G: \mathbb{K}^{n+1} \to \mathbb{K}$, the leading exponent simplifies to $\chi = k$, whereas in case of $k = 1$, we obtain $\chi = dim\, N_2^c$.

Exemplarily, the constellation is qualitatively depicted in figure 4 left with respect to a cusp in $B = \mathbb{R}^3$ and $\bar{z}_l \neq 0, l\ even$. The diagram to the right shows a constellation in $B = \mathbb{R}^3$ with $\bar{z}_l \neq 0, l\ odd$ and the solution curve shifted to one of the coordinate axis, then agreeing with $N_{k+1}$.

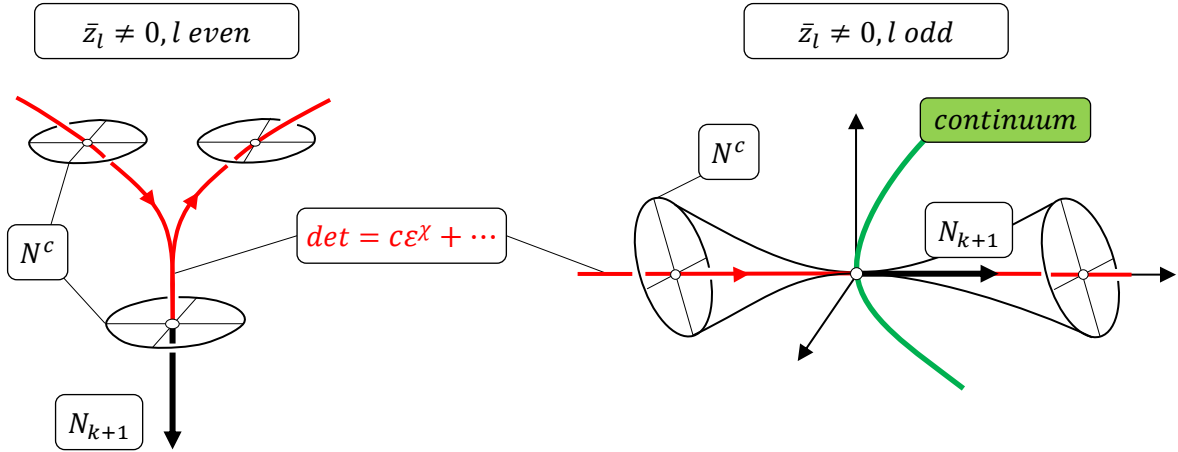

Figure 4 : Two constellations with $B = N^c \oplus N_{k+1}$ and $det$ calculated by the multiplicity $\chi$.

In [9], [18] and [19], the constellation of a given trivial solution curve, as depicted in the right diagram, is treated in detail by use of a direct sum condition of order $k$ comparable to (27). Moreover, the leading exponent $\chi$ (algebraic multiplicity of $k$-transversal eigenvalue) of the determinant is introduced. The theory developed in [9], [10] and [18], [19] represents a powerful generalization of [7] and [20].

Let us now look at an example, how to exploit formulas (29) and (30) with respect to secondary bifurcation based on topological arguments from global bifurcation theory. In particular, decomposition (27) seems to be an extension of the direct sum condition from [9].

For this purpose, consider the equation

$$G[z] = G[\, x, y_1, y_2, y_3\, ] = \begin{pmatrix} y_3 x \\ y_2 + y_1 x \\ y_2 x \end{pmatrix} + H[x,y] = 0 \quad with \quad G: \mathbb{R} \times \mathbb{R}^3 \to \mathbb{R}^3 \quad (31)$$

and $H[x,y] = O(|y|^2)$. For $y = 0$ we obtain a trivial solution curve along the $x$-axis, i.e. in $\varepsilon$-terminology the solution curve simply reads $z(\varepsilon) = \varepsilon \cdot \bar{z}_1 = \varepsilon \cdot (1,0,0,0)^T$, $\bar{z}_i = (0,0,0,0)^T, i \geq 2$. Then, using the iteration from figure 1, we obtain by direct calculation the decomposition

$$B = \mathbb{R}^4 = \left\{\overbrace{\begin{pmatrix}0\\0\\1\\0\end{pmatrix}}^{=N_1^c}\right\} \oplus \left\{\overbrace{\begin{pmatrix}0\\0\\0\\1\end{pmatrix}}^{=N_2^c}\right\} \oplus \left\{\overbrace{\begin{pmatrix}0\\1\\0\\0\end{pmatrix}}^{=N_3^c}\right\} \oplus \left\{\overbrace{\begin{pmatrix}1\\0\\0\\0\end{pmatrix}}^{=N_3}\right\}$$

$$\begin{array}{ccc} \uparrow & \uparrow & \uparrow \\ \boxed{S_1 \neq 0} & \boxed{S_2 \neq 0} & \boxed{S_3 \neq 0} \\ \downarrow & \downarrow & \downarrow \end{array}$$

$$\bar{B} = \mathbb{R}^3 = \underbrace{\left\{\begin{pmatrix}0\\1\\0\end{pmatrix}\right\}}_{=R_1} \oplus \underbrace{\left\{\begin{pmatrix}1\\0\\0\end{pmatrix}\right\}}_{=R_2} \oplus \underbrace{\left\{\begin{pmatrix}0\\0\\1\end{pmatrix}\right\}}_{=R_3} \qquad (32)$$

showing that the trivial solution curve is $k$-regular with degree $k = 2$ and, as usual, with corresponding properties of uniqueness in a pointed wedge $W_0$ combined with stability of $[\bar{z}_1, \bar{z}_2]$ with respect to perturbations of order $O(|x,y|^{2k+1}) = O(|x,y|^5)$. The constellation is qualitatively depicted on the left-hand side within figure 5 below.

In addition, we can now use (29) and (30) to obtain some information about the behavior of the determinant along the $x$-axis. First, the complementary subspace $N^c = N_1^c \oplus N_2^c \oplus N_3^c$ of (32) is obviously given by the $y$-components yielding

$$\chi = 1 \cdot dim\, N_2^c + 2 \cdot dim\, N_3^c = 3 \quad and \quad det\{G_y[\varepsilon,0,0,0]\} = \overset{\neq 0}{\tilde{c}} \cdot \varepsilon^3 + O(|\varepsilon|^4). \qquad (33)$$

Hence, $\chi = 3$ is odd and the determinant with respect to $y$-space, i.e. Brouwer's degree with respect to the pointed wedge $W_0$, changes sign at $\varepsilon = 0$, thus preventing curves of regular values to converge to the $x$-axis. Now, from global bifurcation theory [19], [20], it is well known that a continuum of solutions must emanate from the origin 'absorbing' the limit points of the regular value curves. In the next step, these secondary nontrivial solutions, existing outside the pointed wedge $W_0$, may be looked for by start of another iteration according to figure 1.

In this sense, formulas (29) and (30) may in principle be applied to every solution curve with regularity of degree $k$, possibly indicating by odd leading exponent $l = 1, 3, ...$ of the solution curve in (28), as well as odd leading exponent $\chi = 1, 3, ...$ of the determinant in (30), the existence of further solution curves, perhaps accessible to the iteration from figure 1.

In case of example (31), the nontrivial continuum of solutions depends on higher order terms within perturbation $H[x,y] = O(|y|^2)$, where two possible constellations are shown in figure 5.

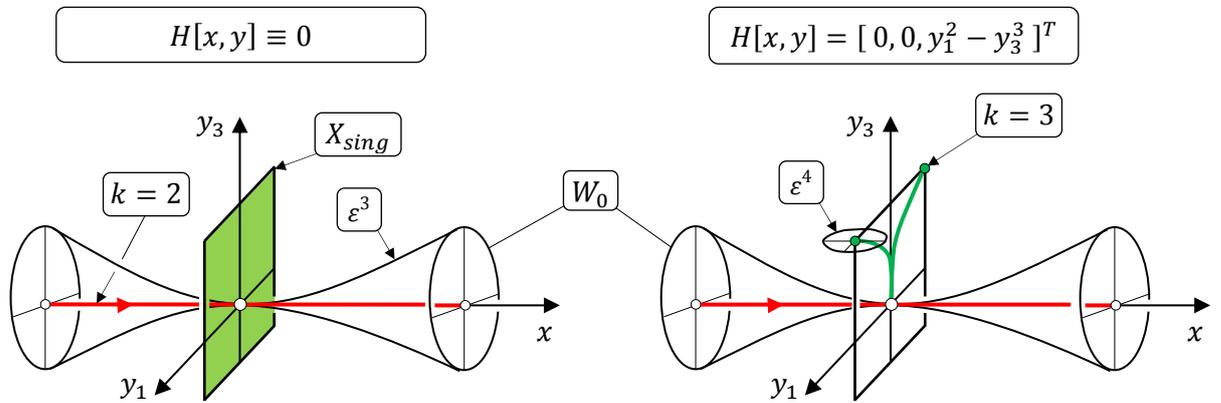

Figure 5 : The nontrivial solutions in case of $H[x,y] \equiv 0$ and $H[x,y] = [\,0, 0, y_1^2 - y_3^3\,]^T$.

In case of $H[x,y] \equiv 0$ (left), the nontrivial continuum of bifurcating solutions is given by the complete $(y_1, y_3)$-plane that is not accessible to the iteration from figure 1 due to instability with respect to perturbations of arbitrary high order. Again, the $(y_1, y_3)$-plane represents the singular locus $X_{sing}$ of (31), consisting of nonisolated critical points with $G'[0, y_1, 0, y_3] = 0$.

In case of $[x, y] = [0, 0, y_1^2 - y_3^3]^T$ (right), a cusp with $k$-regularity of degree $k = 3$ arises within the $(y_1, y_3)$-plane, as can be easily seen by the iteration from figure 1. The pointed wedge of uniqueness of order $O(|\varepsilon|^4)$ belonging to the green marked cusp, is only indicated by an ellipse within figure 5 right-hand side.

For comparison, let us now treat example (31) by use of the direct sum condition from [9] together with the extension from [10] and [18]. Adopting the notation from [9], we obtain

$$G[x, y_1, y_2, y_3] = \overbrace{[G_y^0 + G_{yx}^0 \cdot x]}^{=:L(x)} \cdot y + H[x, y]$$

$$= \left[ \underbrace{\begin{pmatrix} 0 & 0 & 0 \\ 0 & 1 & 0 \\ 0 & 0 & 0 \end{pmatrix}}_{=:L_0} + \underbrace{\begin{pmatrix} 0 & 0 & 1 \\ 1 & 0 & 0 \\ 0 & 1 & 0 \end{pmatrix}}_{=:L_1} \cdot x \right] \cdot y + H[x, y] \qquad (34)$$

and

$$R[L_0] = \left\{ \begin{pmatrix} 0 \\ 1 \\ 0 \end{pmatrix} \right\}, \quad N[L_0] = \left\{ \begin{pmatrix} 1 \\ 0 \\ 0 \end{pmatrix}, \begin{pmatrix} 0 \\ 0 \\ 1 \end{pmatrix} \right\}, \quad L_1 N[L_0] = \left\{ \begin{pmatrix} 0 \\ 1 \\ 0 \end{pmatrix}, \begin{pmatrix} 1 \\ 0 \\ 0 \end{pmatrix} \right\}, \qquad (35)$$

where the direct sum condition from [9] simplifies in case of (31), (34) to the requirement

$$R[L_0] \oplus L_1 N[L_0] = \mathbb{R}^3, \qquad (36)$$

obviously not satisfied by (35). However, in [10] and [18] it is shown that a family of polynomial isomorphisms $\varphi(x)$, $\varphi(0) = I$ can effectively be constructed in such a way that the new family $\bar{L}(x) := L(x)\varphi(x)$ satisfies (36), thus implying the existence of nontrivial solutions as desired. This two stage process can be performed whenever $L_0$ is a Fredholm operator of index zero and 0 is an algebraic eigenvalue of the family $L(x)$ at $x = 0$ [18]. Moreover, in most concrete situations, the second step within the two stage process is not necessary.

Now, the iteration from figure 1 may omit the two stage process, thus generating in an automatic way a direct sum condition that allows for computing solutions of $G[z] = 0$ with respect to existence, uniqueness and stability, as well as the multiplicity $\chi$, possibly implying secondary bifurcation by topological arguments. To some extent, it seems that some essential properties of the isomorphisms from [10], [18] are inherently present within the iteration from figure 1. This aspect proves to be rather convenient when dealing with Banach space lifting in a general context.

The assumption $N_1^c \oplus N_2^c \oplus \cdots \oplus N_{k+1}^c$ to be finite dimensional is not optimal. It is enough to require $N_2^c \oplus \cdots \oplus N_{k+1}^c$ to be finite dimensional and performing a Lyapunov-Schmidt reduction with respect to $G_0^1$ and $N_1^c$, $R_1$ in advance. Note also that the calculation of $\chi$ in (30) does not depend on the dimension $dim\, N_1^c$. If $G_y^0 = L_0$ is a Fredholm operator of index 0, then the subspace $N_2^c \oplus \cdots \oplus N_{k+1}^c$ is of finite dimension.

**Example 4 – Whitney Umbrella :** In this example we show that the iteration is not restricted to deal with isolated solution curves, but can also be used to prove the existence and stability of smooth solution surfaces possessing certain $k$-degree of regularity.

Given the Whitney umbrella, as depicted in figure 6, by the solutions of

$$W[x,y,z] = x^2 - y^2 z = 0 \quad with \quad W : \mathbb{R}^3 \to \mathbb{R} \tag{37}$$

and corresponding $(u,w)$-parametrization of the solution surface

$$x = uw, \quad y = u, \quad z = w^2 \tag{38}$$

complemented by the singular locus $X_{sing}$ along the $z$-axis satisfying $W'[0,0,z] = 0$.

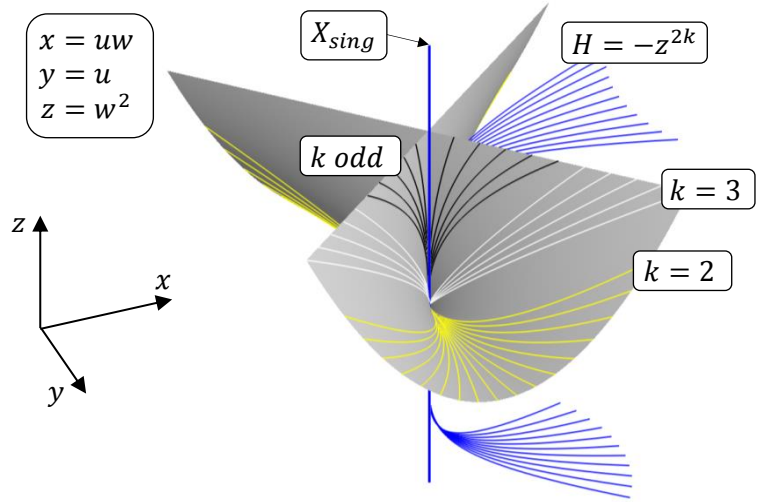

Figure 6 : Whitney umbrella with segments of different $k$-degree.

Again the solutions $X_{sing}$ can be destroyed by an arbitrary small perturbation implying $k$-degree of the $z$-axis to be not finite. On the other hand, there exist different segments within surface (38), showing finite $k$-degree. In some more detail and with respect to perturbations of the Whitney umbrella by $H[x,y,z]$ according to

$$G[x,y,z] = x^2 - y^2 z + H[x,y,z] = 0, \tag{39}$$

the following properties (i)-(iii) follow straightforward from the iteration of figure 1 by use of a parameter dependent version of the implicit function theorem.

(i) If $H[x,y,z] = O(|x,y,z|^3)$ with $H^0_{y^3} = H^0_{y^4} = H^0_{xy^2} = H^0_{zy^2} = 0$, then a smooth $(\varepsilon, v)$-dependent surface of solutions with degree $k = 2$ of regularity of the form

$$\begin{pmatrix} x \\ y \\ z \end{pmatrix} = \varepsilon \cdot \begin{pmatrix} 0 \\ 1 \\ 0 \end{pmatrix} + \varepsilon^2 \cdot \begin{pmatrix} v \\ 0 \\ v^2 \end{pmatrix} + \varepsilon^3 \cdot a(\varepsilon, v) \cdot \begin{pmatrix} 2v \\ 0 \\ -1 \end{pmatrix}$$

exists, comprising a smooth function $a(\varepsilon, v)$ with $v$ bounded and $\varepsilon$ chosen sufficiently small. In case of $H \equiv 0$ (pure Whitney umbrella), we obtain $a(\varepsilon, v) \equiv 0$ yielding the yellow curves within figure 6 that are tangentially touching the $y$-axis in the origin. Under perturbation by $H \neq 0$, these curves are slightly deformed in $x$- and $z$-direction by terms of order $O(|\varepsilon|^3)$. Within the iteration of figure 1, surjectivity of the sum operator is reached in one step by $S_3 = 6G_0^3 \bar{z}_1^2 + 6G_0^2 \bar{z}_2 = 6 [4v \; 0 \; -2]$.

Note also that in general, a solution curve with $k = 2$ only is stable with respect to perturbations of order $2k + 1 = 5$. The improvement to most of the monomials of order 3 and order 4 results from some simplifying structure inherent within defining equation (37) of the Whitney umbrella.

(ii) If $H[x,y,z] = O(|x,y,z|^{2k+1})$ and $k = 3, 5, ..., odd$, then a smooth $(\varepsilon, v)$-dependent surface of solutions with degree $k$ of regularity of the form

$$\begin{pmatrix} x \\ y \\ z \end{pmatrix} = \varepsilon^2 \cdot \begin{pmatrix} 0 \\ 0 \\ 1 \end{pmatrix} + \varepsilon^{k-1} \cdot \begin{pmatrix} 0 \\ v \\ 0 \end{pmatrix} + \varepsilon^k \cdot \begin{pmatrix} v \\ 0 \\ 0 \end{pmatrix} + \varepsilon^{k+1} \cdot a(\varepsilon, v) \cdot \begin{pmatrix} 1 \\ 0 \\ 0 \end{pmatrix}$$

exists, comprising a smooth function $a(\varepsilon, v)$ with $v \neq 0$ bounded and $\varepsilon$ chosen sufficiently small. Again, in case of pure Whitney umbrella $H \equiv 0$, we obtain $a(\varepsilon, v) \equiv 0$ implying white curves in figure 6 with $k = 3$ and black curves with $k \geq 5$, all touching the $z$-axis in the origin by order of $(k-3)/2$. Surjectivity of the sum operator is reached by $S_{k+1} = c \cdot G_0^2 \bar{z}_k = c \cdot [2v \ 0 \ 0]$, $c \neq 0$, $v \neq 0$.

(iii) If $H[x,y,z] = z^{2k} \cdot r(x,y,z)$, $r(0,0,0) < 0$ and $k = 1, 2, 3, ...$, then two $(\varepsilon, v)$-dependent surfaces of solutions of the form

$$\begin{pmatrix} x \\ y \\ z \end{pmatrix} = \varepsilon \cdot \begin{pmatrix} 0 \\ 0 \\ 1 \end{pmatrix} + \varepsilon^k \cdot \begin{pmatrix} \pm 1 \\ v \\ 0 \end{pmatrix} + \varepsilon^{k+1} \cdot a_\pm(\varepsilon, v) \cdot \begin{pmatrix} 1 \\ 0 \\ 0 \end{pmatrix}$$

exist with degree $k$ of regularity and smooth functions $a_\pm(\varepsilon, v)$ with $v$ bounded and $\varepsilon$ sufficiently small. Surjectivity of the sum operator is reached by $S_{k+1} = c \cdot G_0^2 \bar{z}_k = c \cdot [2 \ 0 \ 0]$, $c \neq 0$. These surfaces arise from the $z$-axis $X_{sing}$ of the Whitney umbrella under perturbation by $H \leq 0$ of order $O(|z|^{2k})$, touching the $z$-axis by order of $k - 1$. In figure 6, only one of the two surfaces is indicated by blue lines with $k = 8$.

Note also that in case of $r(x,y,z) \equiv -1$, equation (39) turns into $G[x,y,z] = x^2 - y^2 z - z^{2k} = 0$ yielding by splitting-lemma [1], [4] the normal form of a $D_{2k+1}$-singularity in $\mathbb{R}^3$.

**Example 5 – G. Belitskii and D. Kerner :** Consider example 4.1 from [2]

$$G[x_1, x_2, y_1, y_2] = H[x_1, x_2, y_1, y_2] + y_1 x_1^k + y_2 x_2^k + p[x_1, x_2] = 0 \quad (40)$$

under the restriction $G: \mathbb{R}^4 \to \mathbb{R}$ and $H[x, y] = O(|y|^2)$, $k > 2$. Then, in case of $p[x_1, x_2] = 0$, equation (40) implies basic solutions $y \equiv 0$ that may be destroyed when choosing $p[x_1, x_2] \neq 0$.

Now in [2], it is shown by application of the implicit function theorems of Tougeron [30] and Fisher [11] that the basic solutions $y \equiv 0$ are continued to a smooth function $y(x)$ with $y(0) = 0$, as long as the perturbation $p[x_1, x_2]$ of order $2k + 1$ is restricted to

$$p[x_1, x_2] = O(|x_1|^{2k+1} + |x_1|^{k+1}|x_2|^k + |x_1||x_2|^{2k} + |x_2||x_1|^{2k} + |x_1|^k|x_2|^{k+1} + |x_2|^{2k+1}).$$

The result is further improved by the strong implicit function theorem for modules, established in [2], to arbitrary perturbations $p[x_1, x_2]$ of order $O(|x_1, x_2|^{2k+1})$.

Now, when applying the iteration from figure 1 to equation (40), then the following result is obtained with respect to $(x, y)$-dependent perturbations of the form

$$p[x_1, x_2, y_1, y_2] = O(|x_1, x_2, y_1, y_2|^{2k+1}). \quad (41)$$

(i) There exists a smooth 3-dimensional surface of solutions

$$\begin{pmatrix} x_1 \\ x_2 \\ y_1 \\ y_2 \end{pmatrix} = \varepsilon \cdot \begin{pmatrix} \cos\varphi \\ \sin\varphi \\ 0 \\ 0 \end{pmatrix} + \varepsilon^{k+1} \cdot \left[ a(\varepsilon, \varphi, b) \cdot \begin{pmatrix} 0 \\ 0 \\ \cos^k \varphi \\ \sin^k \varphi \end{pmatrix} + b \cdot \begin{pmatrix} 0 \\ 0 \\ -\sin^k \varphi \\ \cos^k \varphi \end{pmatrix} \right] \quad (42)$$

with degree $k$ of regularity and a smooth function $a(\varepsilon, \varphi, b)$, $\varphi \in [0, 2\pi[$, $b$ bounded and $\varepsilon$ chosen sufficiently small. Surjectivity of the sum operator is reached in one step by $S_{k+1} = c \cdot G_0^{k+1} \bar{z}_1^k = c \cdot [0\ 0\ cos^k\varphi\ sin^k\varphi]$, $c \neq 0$.

To omit multiple parametrization, we further restrict to $\varepsilon \geq 0$. Note that perturbations by (41) are allowed to include terms of order $O(|x_i|^{2k}|y_j|)$, $i, j = 1, 2$ that may not be allowed in [2]. Terms of $p[x, y]$, quadratic in $y$, may be absorbed by $H[x, y] = O(|y|^2)$ in (40).

In case of $p[x, y] = 0$, the function $a(\varepsilon, \varphi, b)$ satisfies $a(\varepsilon, \varphi, 0) = 0$, hence yielding the basic solution $y(\varepsilon, \varphi, 0) \equiv 0$ with $b = 0$ and the $(x_1, x_2)$-plane parametrized by polar coordinates according to $(x_1, x_2) = \varepsilon(cos\ \varphi, sin\ \varphi)$. The basic solution is indicated in the left diagram of figure 7 by the red plane and labeled with $b = 0$. When choosing $b \neq 0$, the basic solution is modified by terms of order $O(|\varepsilon|^{k+1})$ according to (42), as depicted by grey and yellow surfaces in the left diagram.

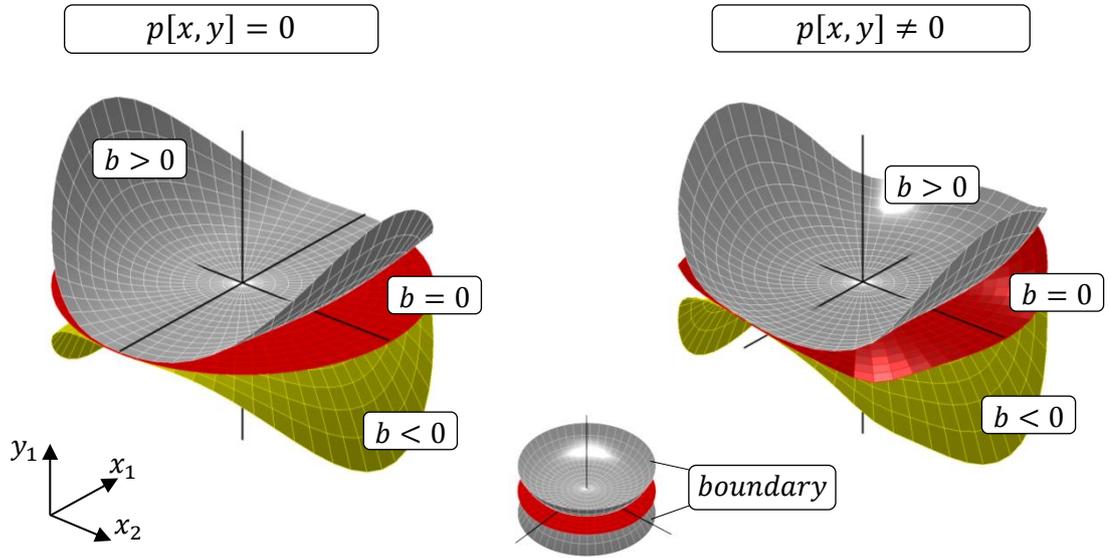

Figure 7 : The solution surfaces in case of $p[x, y] = 0$ (left) and $p[x, y] \neq 0$ (right).

Under perturbation by $p[x, y] \neq 0$ and according to (41), the solution surfaces of the left diagram are typically modified as qualitatively indicated in the right diagram. In this sense, the solutions of the left diagram are continued under perturbations of the form (41). In particular, all of the surfaces in the right diagram are passing through the origin by choice of $\varepsilon = 0$ in (42), thus repeating the result from [2].

Additionally, by uniqueness of the implicit function theorem, the 3-dimensional surface (42) represents all solutions of equation (40) within an open neighborhood of $(x_1, x_2)$-plane in $\mathbb{R}^4$, defined by (42) with $a(\varepsilon, \varphi, b)$ replaced by a parameter $a$ varying in a bounded interval, i.e. the neighborhood is defined by four parameters $(\varepsilon, \varphi, a, b)$, as qualitatively indicated in the middle diagram of figure 7. Here, grey surfaces indicate upper and lower boundary of the open neighborhood of uniqueness within $(x_1, x_2, y_1)$-space. The two external parameters $(\varepsilon, \varphi)$ in (42) may be substituted by an internal $(x_1, x_2)$-parametrization.

*3. Results*

First we restrict to pure existence, giving an Artin-Tougeron type result in Banach spaces.

As preliminary, we summarize the state of the iteration at the red dot in figure 1. The leading coefficients $[\bar{z}_1, \ldots, \bar{z}_k] \in \pi_k(\bar{X}_{2k})$, the affine subspace $\bar{X}_{2k}$ by itself, and the sequence $B = N_0 \supset N_1 \supset \cdots \supset N_k \supset N_{k+1}$, $k \geq 1$ are built up by solvability and decomposition conditions yielding

$$B = N_1^c \oplus N_2^c \oplus \cdots \oplus N_{k+1}^c \oplus N_{k+1}$$
$$\uparrow \quad \uparrow \quad\quad \uparrow$$
$$\boxed{S_1} \quad \boxed{S_2} \quad\quad \boxed{S_{k+1}} \tag{43}$$
$$\downarrow \quad \downarrow \quad\quad \downarrow$$
$$\bar{B} = R_1 \oplus R_2 \oplus \cdots \oplus R_{k+1} \oplus R_{k+1}^c$$

with corresponding sum operator

$$[S_1 \cdots S_{k+1}] \begin{pmatrix} n_0 \\ \vdots \\ n_k \end{pmatrix} = S_1 n_0 + \cdots + S_{k+1} n_k \in L[N_0 \times \cdots \times N_k, \bar{B}] \tag{44}$$

only depending on the leading coefficients $[\bar{z}_1, \ldots, \bar{z}_k]$ and the first $k+1$ derivatives of $G$ at $z=0$ denoted by $G_0^1, \ldots, G_0^{k+1}$.

**Theorem 1 - Existence**

Assume $[S_1 \cdots S_{k+1}] \in L[N_0 \times \cdots \times N_k, \bar{B}]$ to be surjective, i.e. $R_{k+1}^c = \{0\}$.

Then a smooth family $z(\varepsilon, n_1, \ldots, n_{k+1})$, $|\varepsilon| \ll 1$, $n_i \in N_i$, $i = 1, \ldots, k+1$ of solution curves of $G[z] = 0$ exists with the following properties

(i) $z(\varepsilon, \bar{n}) = \varepsilon \bar{z}_1 + \cdots + \frac{1}{k!} \varepsilon^k \bar{z}_k + \frac{1}{(k+1)!} \varepsilon^{k+1} z_{k+1}(\varepsilon, \bar{n}) + \cdots + \frac{1}{(2k+1)!} \varepsilon^{2k+1} z_{2k+1}(\varepsilon, \bar{n})$

with $(\varepsilon, \bar{n}) := (\varepsilon, n_1, \ldots, n_{k+1})$ and $\bar{n}$ bounded in $N_1 \times \cdots \times N_{k+1}$.

(ii) $[\bar{z}_1, \ldots, \bar{z}_k, z_{k+1}(0, \bar{n}), \ldots, z_{2k}(0, \bar{n}), z_{2k+1}(0, \bar{n})] \in \bar{X}_{2k+1}$

$[\bar{z}_1, \ldots, \bar{z}_k, z_{k+1}(\varepsilon, \bar{n}), \ldots, z_{2k}(\varepsilon, \bar{n})] \in \bar{X}_{2k}$

(iii) $\begin{pmatrix} z_{2k+1} \\ \vdots \\ z_{k+1} \end{pmatrix}(\varepsilon, \bar{n}) = \hat{I}_{k+1} + \boxed{\widehat{M}_{k+1}} \cdot \left[ \begin{pmatrix} n_1^c \\ \vdots \\ n_{k+1}^c \end{pmatrix}(\varepsilon, \bar{n}) + \boxed{\hat{L}_{k+1}} \cdot \begin{pmatrix} n_1 \\ \vdots \\ n_{k+1} \end{pmatrix} \right]$

with constant vector $\hat{I}_{k+1} \in B^{k+1}$, constant matrices $\boxed{\widehat{M}_{k+1}}$, $\boxed{\hat{L}_{k+1}} \in GL[B^{k+1}, B^{k+1}]$
and $(n_1^c, \ldots, n_{k+1}^c)(0, \bar{n}) = 0$.

(iv) Assume $N^c := N_1^c \oplus N_2^c \oplus \cdots \oplus N_{k+1}^c$ to be of finite dimension.

Then the determinant of the linearization $G_{n^c}[z(\varepsilon, \bar{n})] \in L[N^c, \bar{B}]$ satisfies

$$det\{G_{n^c}[z(\varepsilon, \bar{n})]\} = c \cdot \varepsilon^\chi + O(|\varepsilon|^{\chi+1}), \quad c \neq 0$$

with leading exponent

$$\chi = 1 \cdot \dim N_2^c + \cdots + k \cdot \dim N_{k+1}^c.$$

The elements within $\hat{I}_{k+1}$, $\widehat{M}_{k+1}$ and $\hat{L}_{k+1}$ are constructively defined by the composition of multilinear mappings beween subspaces of the decomposition (43).

**Remarks 1)** Once the iteration works, the proof of Theorem 1 is straightforward. The affine subspace $\bar{X}_{2k}$ is plugged into the Ansatz

$$G\left[\varepsilon \bar{z}_1 + \cdots + \frac{1}{k!}\varepsilon^k \bar{z}_k + \frac{1}{(k+1)!}\varepsilon^{k+1} z_{k+1} + \cdots + \frac{1}{(2k+1)!}\varepsilon^{2k+1} z_{2k+1}\right] = 0 \quad (45)$$

yielding by construction a remainder equation of the form

$$T^{2k+1}\left[\underbrace{\bar{z}_1, \ldots, \bar{z}_k, z_{k+1}, \ldots, z_{2k}}_{\in \bar{X}_{2k}}, z_{2k+1}\right] + O(|\varepsilon|) = 0. \quad (46)$$

Then with $\varepsilon = 0$, the equation $T^{2k+1}[\bar{X}_{2k}, z_{2k+1}] = 0$ is easily solved by use of surjectivity of the sum operator $[S_1 \cdots S_{k+1}]$, implying the existence of an $\bar{n}$-dependent family of basic solutions (in fact defining $\bar{X}_{2k+1}$ in (ii)) that are locally continued to $\varepsilon \neq 0$ by surjectivity and the implicit function theorem. In this sense, a subspace of $\bar{X}_{2k}$ is first lifted to $\bar{X}_{2k+1}$ and further lifted to $\bar{X}_\infty$ by implicit function theorem. Here $\bar{X}_\infty \subset X_\infty$ denotes the elements in arc space $X_\infty$ with leading coefficients $[\bar{z}_1, \ldots, \bar{z}_k]$.

Note also that by construction, the components $[\bar{z}_1, \ldots, \bar{z}_k, z_{k+1}(\varepsilon, \bar{n}), \ldots, z_{2k}(\varepsilon, \bar{n})]$ remain in $\bar{X}_{2k}$ for all values of $\varepsilon, |\varepsilon| \ll 1$, as stated in the second line of (ii).

The special structure of higher order coefficients $[z_{k+1}, \ldots, z_{2k+1}]$ given in (iii) is adopted from a filtration within $\bar{X}_{2k}$ and corresponding decompositions of $B$ and $\bar{B}$ according to (43). In some more detail, the family of solutions in (i) is parametrized by the kernels of $S_1, \ldots, S_{k+1}$

$$N_1 \supset \cdots \supset N_k \supset N_{k+1} \quad (47)$$

with dependent variables chosen from corresponding subspaces

$$N_1^c \oplus N_2^c \oplus \cdots \oplus N_{k+1}^c = B/N_1 \oplus N_1/N_2 \oplus \cdots \oplus N_k/N_{k+1}. \quad (48)$$

Note that the parametrization of the solution curves obviously depends on the choice of the subspaces $N_i^c$, where it can be shown that $k$-regularity is an invariant with respect to chosen subspaces. In addition, the quotient spaces in (48) indicate that Theorem 1 may also be formulated as a coordinate free version.

Finally, the determinant in (iv) offers the possibility to calculate Brouwer's degree along the solution curves, thus allowing for search of global bifurcation phenomena [19], [20], as exemplified within example 3 of section 2.

Theorem 1 may be reformulated with respect to differentiable $C^n$-maps $G: B \to \bar{B}$ between Banach spaces with $n$ chosen sufficiently high.

The structure of the proof is given in [24], details can be found in [25]. In [24], the focus is layed on singularly perturbed nonhyperbolic points and applications to biomathematics in molecular cell biology [29], whereas [25] concentrates on mathematical aspects of bifurcation theory based on the system of undetermined coefficients.

**2)** It should be noted that instead of (46) we can also work with the remainder equation

$$T^{2k}[\bar{z}_1, \ldots, \bar{z}_{k-1}, z_k, \ldots, z_{2k}] + O(|\varepsilon|) = 0 \quad (49)$$

implying solutions of the form

$$G\left[ \varepsilon \bar{z}_1 + \cdots + \frac{1}{(k-1)!}\varepsilon^{k-1}\bar{z}_{k-1} + \frac{1}{k!}\varepsilon^k z_k(\varepsilon,\bar{n}) + \cdots + \frac{1}{(2k)!}\varepsilon^{2k} z_{2k}(\varepsilon,\bar{n}) \right] = 0 \qquad (50)$$

which represent a certain generalization of Theorem 1, because the coefficient $z_k = z_k(\varepsilon,\bar{n})$ is now allowed to vary in the vicinity of $\bar{z}_k$. But also when performing this extension, first a basic solution is constructed in $\pi_{2k}(\bar{X}_{2k+1})$ and further continued to $\bar{X}_\infty$ using surjectivity of the same operator $[S_1 \cdots S_{k+1}]$. We preferred (46) because of possible reference to [24], [25].

**3)** When reaching the red dot in figure 1, the affine subspace $\bar{X}_{2k} \subset B^{2k}$ is successfully constructed, yielding an infinity of approximate solution curves satisfying $G[z] = 0$ up to order $O(|\varepsilon|^{2k+1})$, i.e. every $(2k)$-tuple $[\bar{z}_1, \ldots, \bar{z}_k, z_{k+1}, \ldots, z_{2k}] \in \bar{X}_{2k}$ gives rise to an approximate solution curve such that

$$G\left[ \varepsilon\bar{z}_1 + \cdots + \frac{1}{k!}\varepsilon^k \bar{z}_k + \frac{1}{(k+1)!}\varepsilon^{k+1} z_{k+1} + \cdots + \frac{1}{(2k)!}\varepsilon^{2k} z_{2k} \right] = O(|\varepsilon|^{2k+1}). \qquad (51)$$

In general, these approximate solution curves only agree with the exact solution curves $z(\varepsilon,\bar{n})$ from (i) with respect to first $k$ leading coefficients $[\bar{z}_1, \ldots, \bar{z}_k]$. This is well-known from Artin-Approximation that the approximate solution curve loses some derivatives compared to the lifted, exact solutions. In our context, the last $k$ derivatives are in general lost, whereas the first $k$ derivatives are maintained. In some more detail, this aspect is treated in Corollary 2 below.

**4)** It is interesting to perform the iteration in the simplifying case of $G[z]$ to be a polynomial. One may hope that the solvability conditions in figure 1 turn out to be satisfied automatically when $k$ exceeds some specific value, possibly implying a result comparable to Strong Artin- or Greenberg-Approximation [22]. In fact some recurrency structure is quite obvious, but we did not succeed to eliminate the obstruction given by the solvability conditions. Some questions concerning weak Greenberg functions and lower bounds of Greenberg functions are treated in the remarks of Corollary 2 below.

**5)** Consider the family of solution curves $z(\varepsilon,\bar{n})$ from (i) with $k$-regularity of $[\bar{z}_1, \cdots, \bar{z}_k]$. Then $k$-regularity of leading coefficients should be invariant with respect to diffeomorphic $\varepsilon$-parameter transformations as well as diffeomorphic $z$-coordinate transformations. In some more detail, the left/right-transformation of the solution curves

$$\hat{z}(\varepsilon,\bar{n}) := \varphi^{-1} \circ z(\,\cdot\,,\bar{n}) \circ \psi(\varepsilon) = \varphi^{-1}[\,z(\psi(\varepsilon),\bar{n})\,] \qquad (52)$$

with $\varphi : B \to B$, $\varphi(0) = 0$ diffeomorphic near $z = 0$ and $\psi : \mathbb{K} \to \mathbb{K}$, $\psi(0) = 0$ diffeomorphic near $\varepsilon = 0$ should imply solution curves of the transformed map $\hat{G}[z] := G[\,\varphi(z)\,]$ with same $k$-degree of regularity. Note also that the $k$-tuple of leading coefficients $[\hat{z}_1, \cdots, \hat{z}_k]$ of the transformed solution curves $\hat{z}(\varepsilon,\bar{n})$ can easily be computed from $[\bar{z}_1, \cdots, \bar{z}_k]$ using derivatives of left/right transformation at $z = 0$ and $\varepsilon = 0$ respectively.

Now, a first analysis strongly suggests this invariance of $k$-regularity with respect to left/right-transformation to hold true. At least, if we restrict to transformations of the form $\varphi(z) = z + O(|z|^{k+1})$ and $\psi(\varepsilon) = \varepsilon + O(|\varepsilon|^{k+1})$, i.e. if the $k$-jets of both transformations equal identity, then the leading coefficients $[\bar{z}_1, \cdots, \bar{z}_k]$ remain unchanged and invariance of $k$-regularity is easily established.

**6)** If the Banach space $B$ splits according to $B = \mathbb{K} \times Y$ and $Y$ Banach space with respect to $\mathbb{K}$, then the solution curves may be splitted too by

$$z(\varepsilon,\bar{n}) = \begin{pmatrix} x(\varepsilon,\bar{n}) \\ y(\varepsilon,\bar{n}) \end{pmatrix} = \begin{pmatrix} \frac{1}{l!}\varepsilon^l \cdot \bar{x}_l + O(|\varepsilon|^{l+1}) \\ \frac{1}{\tau!}\varepsilon^\tau \cdot \bar{y}_\tau + O(|\varepsilon|^{\tau+1}) \end{pmatrix}$$

with $\bar{x}_l, \bar{y}_\tau \neq 0$, $l, \tau \geq 1$ denoting the leading coefficients of $x$- and $y$-expansion respectively. An $\varepsilon$-transformation of the form $\psi(\varepsilon, \bar{n}) = \varepsilon + \varepsilon^2 \cdot \bar{\psi}(\varepsilon, \bar{n})$ can be performed to eliminate higher order terms $O(|\varepsilon|^{l+1})$ of the $x$-expansion, yielding new parametrization such that

$$\begin{pmatrix} x(\psi(\varepsilon,\bar{n}),\bar{n}) \\ y(\psi(\varepsilon,\bar{n}),\bar{n}) \end{pmatrix} = \begin{pmatrix} \frac{1}{l!}\varepsilon^l \cdot \bar{x}_l \\ \frac{1}{\tau!}\varepsilon^\tau \cdot \bar{y}_\tau + O(|\varepsilon|^{\tau+1}) \end{pmatrix}.$$

This kind of Puiseux parametrization of the solution curves by the external parameter $\varepsilon$ may further be brought to an internal $x$-parametrization by solving $x = \frac{1}{l!}\varepsilon^l \cdot \bar{x}_l$ with respect to $\varepsilon$, finally implying the $y$-component to be parametrized internally by $x$ according to

$$y(x,\bar{n}) = x^{\frac{\tau}{l}} \cdot \tilde{y}_\tau + O\left(|x|^{\frac{\tau+1}{l}}\right)$$

with common denominator $l \geq 1$ of all exponents. This kind of internal parametrization is directly obtained when working with Newton-polygons, as investigated in some more detail within section 4. Note also, if $l = 1$, then the $y$-coordinate transformation $y = y(x, \bar{n}) + \bar{y}$ implies

$$\bar{G}[x, \bar{y}, \bar{n}] := G[x, y(x, \bar{n}) + \bar{y}] = 0 \quad \text{with} \quad |x| \ll 1$$

and we obtain trivial solutions given by $\bar{G}[x, 0, \bar{n}] = 0$, i.e. the solutions are transformed to the $x$-axis, typically representing the starting constellation of global bifurcation theory addressed within example 3 of section 2.

**7)** Another possibility to restrict to internal $x$-parametrization in case of $B = \mathbb{K} \times Y$ is given by restricting the iteration of figure 1 to an ansatz of the form

$$z = \varepsilon \cdot \underbrace{\begin{pmatrix} 1 \\ y_1 \end{pmatrix}}_{=z_1} + \frac{1}{2}\varepsilon^2 \cdot \underbrace{\begin{pmatrix} 0 \\ y_2 \end{pmatrix}}_{=z_2} + \frac{1}{6}\varepsilon^3 \cdot \underbrace{\begin{pmatrix} 0 \\ y_3 \end{pmatrix}}_{=z_3} + \cdots$$

characterized by $z_1 = (1, y_1)^T$ and $z_i = (0, y_i)^T$, $i \geq 2$, yielding solution curves of the form $z(\varepsilon) = [\varepsilon, y(\varepsilon)]$ or $z(x) = [x, y(x)]$ respectively. The corresponding system of undetermined coefficients

$$T^1\left[\begin{pmatrix} 1 \\ y_1 \end{pmatrix}\right] = 0$$
$$\vdots$$
$$T^n\left[\begin{pmatrix} 1 \\ y_1 \end{pmatrix}, \ldots, \begin{pmatrix} 0 \\ y_n \end{pmatrix}\right] = 0$$

defines as $n \to \infty$ some sort of reduced arc space $X_\infty^0 \subset X_\infty$ with formal approximations that may correspond to $x$-parametrizations of exact solution curves. Then, using the iteration from figure 1, $k$-regularity of $[\bar{z}_1, \cdots, \bar{z}_k] = [(1, \bar{y}_1)^T, (0, \bar{y}_2)^T, \ldots, (0, \bar{y}_k)^T]$ can be proven by an adapted decomposition of $B$ according to

$$B = \mathbb{K} \times Y = \overbrace{\begin{pmatrix} 0 \\ Y_1^c \end{pmatrix}}^{=N_1^c} \oplus \cdots \oplus \overbrace{\begin{pmatrix} 0 \\ Y_{k+1}^c \end{pmatrix}}^{=N_{k+1}^c} \oplus \overbrace{\begin{pmatrix} 0 \\ Y_{k+1} \end{pmatrix}}^{=N_{k+1}} \oplus \left\{\begin{pmatrix} 1 \\ \bar{y}_1 \end{pmatrix}\right\} \tag{53}$$

with appropriate subspaces $Y_i^c$, $i = 1, \ldots, k+1$ and $Y_{k+1}$ of $Y$, i.e. a specific version of Theorem 1 arises. The corresponding solution sets of the system of undetermined coefficients, characterized by fixed leading coefficients $[(1, \bar{y}_1)^T, (0, \bar{y}_2)^T \ldots, (0, \bar{y}_k)^T]$ and $z_i = (0, y_i)^T$, $i \geq k+1$, are further labelled by $\bar{X}_{2k}^0$, $\bar{X}_{2k+1}^0$ and $\bar{X}_\infty^0$. These sets are investigated more closely after Corollary 2 below.

Next we turn to some questions of perturbation and stability. At the red dot in figure 1, the iteration process has dealt with equations $T^1 = \cdots = T^{2k} = 0$ depending on the derivatives $G_0^1, \ldots, G_0^{2k}$. This indicates that the solution curves of Theorem 1 should be stable with respect to perturbations of $G[z]$ of order $O(|z|^{2k+1})$.

Conversely, we may conclude that solution curves which can be destroyed by perturbations of arbitrary high order (or even by a flat map) cannot be handled by the iteration. But apart from this situation, typically characterized by a singular locus of nonisolated critical points, we hope to grasp all branches of finite type by use of the iteration in figure 1.

Now, if we look in detail on the iteration process, the following properties with respect to perturbation and stability can easily be summarized.

**Corollary 1 - Stability**

(i) Perturbations of the form $G[z] + H[z] = 0$ with $H[z] = O(|z|^{2k+1})$ leave Theorem 1 unchanged, only varying derivatives of higher order coefficients $z_{k+1}(\varepsilon, \bar{n}), \ldots, z_{2k+1}(\varepsilon, \bar{n})$.

(ii) Perturbations of the form $G[z] + c \cdot H[z] = 0$ with $H[z] = O(|z|^{2k})$ do not destroy the family of solution curves if the constant $c \in \mathbb{K}$ is chosen sufficiently small.

(iii) Perturbations of the form $G[z] + H[z] = 0$ with $H[z] = O(|z|^l)$, $l \leq 2k - 1$ may destroy the solution curves of Theorem 1.

**Remarks 1)** As already mentioned, the sum operator (44) only depends on the derivatives $G_0^1, \ldots, G_0^{k+1}$, implying independence of surjectivity with respect to perturbations of $G[z]$ of order $O(|z|^{k+2})$. Hence the reason for destruction of solution branches by perturbations of order $O(|z|^l)$, $l \leq 2k - 1$ in (iii) is not caused by loss of surjectivity, but by the solvability condition that cannot be satisfied anymore.

**2)** Let us denote by $J_l, l \geq 0$ the set of maps agreeing with $G[z]$ up to order $l$, i.e. the $l$-jets of maps within $J_l$ are supposed to be identical. Further, let us denote by $R_l, l \geq 0$ the set of maps to be right equivalent to $G[z]$ with respect to coordinate transformations $\varphi(z)$ with $l$-jet the identity, i.e. $\varphi(z) = z + O(|z|^{l+1})$.

Then, according to Corollary 1 (i), $k$-regularity of $[\bar{z}_1, \cdots, \bar{z}_k]$ is an invariant of $J_{2k}$, as well as an invariant of $R_k$, as expected according to Remark 5) of Theorem 1. In addition, the relations $R_k \subset R_0$, $R_k \subset J_k$ and $J_{2k} \subset J_k$ are obviously valid, implying a basic constellation as qualitatively indicated in the upper diagram of figure 8.

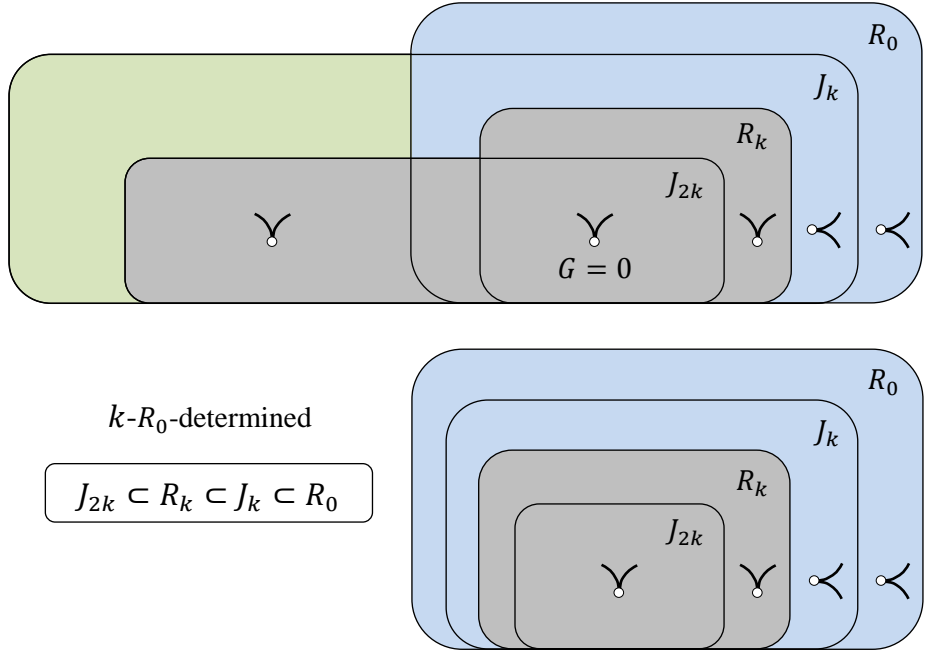

Figure 8 : Behavior of solution curve within different sets of maps.

Note that a solution curve with $k$-regularity of $[\bar{z}_1, \cdots, \bar{z}_k]$ is maintained within the grey set $J_{2k} \cup R_k$, as exemplarily indicated by a cusp curve. Within the blue set $R_0$, the solution curve typically occurs with transformed leading coefficients $[\hat{z}_1, \cdots, \hat{z}_k]$ and $k$-regularity expected again by Remark 5) of Theorem 1. Within the green set $J_k \setminus (R_0 \cup J_{2k})$, the solution curve is typically destroyed.

Now, if we additionally assume, $G : \mathbb{R}^n \to \mathbb{R}$ to be a real function, as well as $J_k \subset R_0$, i.e. $G$ is supposed to be $k$-$R_0$-determined, then $G$ follows to be $(2k)$-$R_k$-determined [5], ensuring the simple sequence $J_{2k} \subset R_k \subset J_k \subset R_0$ to hold true. In this sense, the strong property of $G$ to be $k$-$R_0$-determined forces both of the sets $J_k$ and $J_{2k}$ to move into $R_0$ and $R_k$ respectively, ending up with a rather simple constellation as depicted in the lower diagram of figure 8. In particular, the green set without solution curve has completely vanished.

In the next step, we aim to complement Theorem 1 with respect to some relations between derivatives of approximate and lifted solution curves, as well as a relation between the sets $\bar{X}_{2k+l}$ and $\bar{X}_\infty$.

As already mentioned in Remark 3 of Theorem 1, an approximate solution curve defined by $[\bar{z}_1, \ldots, \bar{z}_k, z_{k+1}, \ldots, z_{2k}] \in \bar{X}_{2k}$, typically loses the last $k$ derivatives compared to the lifted exact solution curve. Now, when looking at Tougeron's implicit function theorem [22] or Hensel's Lemma [11], one observes that an improvement of the approximation given by $[\bar{z}_1, \ldots, \bar{z}_k, z_{k+1}, \ldots, z_{k+l}, z_{k+l+1}, \ldots, z_{2k+l}] \in \bar{X}_{2k+l}$, $l \geq 1$ and

$$G\left[\varepsilon\bar{z}_1 + \cdots + \frac{1}{k!}\varepsilon^k \bar{z}_k + \frac{1}{(k+1)!}\varepsilon^{k+1} z_{k+1} + \cdots + \frac{1}{(2k+l)!}\varepsilon^{2k+l} z_{2k+l}\right] = O(|\varepsilon|^{2k+l+1}) \quad (54)$$

should give rise to an exact solution curve that agrees in the coefficients $[\bar{z}_1, \ldots, \bar{z}_k, z_{k+1}, \ldots, z_{k+l}]$ with the approximate solution from (54). This means that an improvement of the approximation effectively improves the accordance between both curves. Or to be more precise, only the last $k$ coefficients $[z_{k+l+1}, \ldots, z_{2k+l}]$ have to be changed, when lifting the approximation to an exact solution. In context of proof, only the last $k$ coefficients are needed as variables for use within the implicit function theorem, as can easily be seen along the following lines.

First note that if $[S_1 \cdots S_{k+1}]$ becomes surjective at the red point in figure 1, we can decide whether to apply at once the implicit function theorem for lifting the approximation $[\bar{z}_1, \ldots, \bar{z}_k]$ to $\bar{X}_\infty$ or to continue the iteration and eventually apply the implicit function theorem on a later stage of the process.

In addition, when continuing the iteration, we will be able to compute sequentially the complete sets $\bar{X}_{2k+1}, \bar{X}_{2k+2}, \ldots, \bar{X}_{2k+l}$ with $l$ arbitrary large. This results from the observation that the iteration process simplifies heavily, as soon as $[S_1 \cdots S_{k+1}]$ has become surjective, as indicated in figure 9.

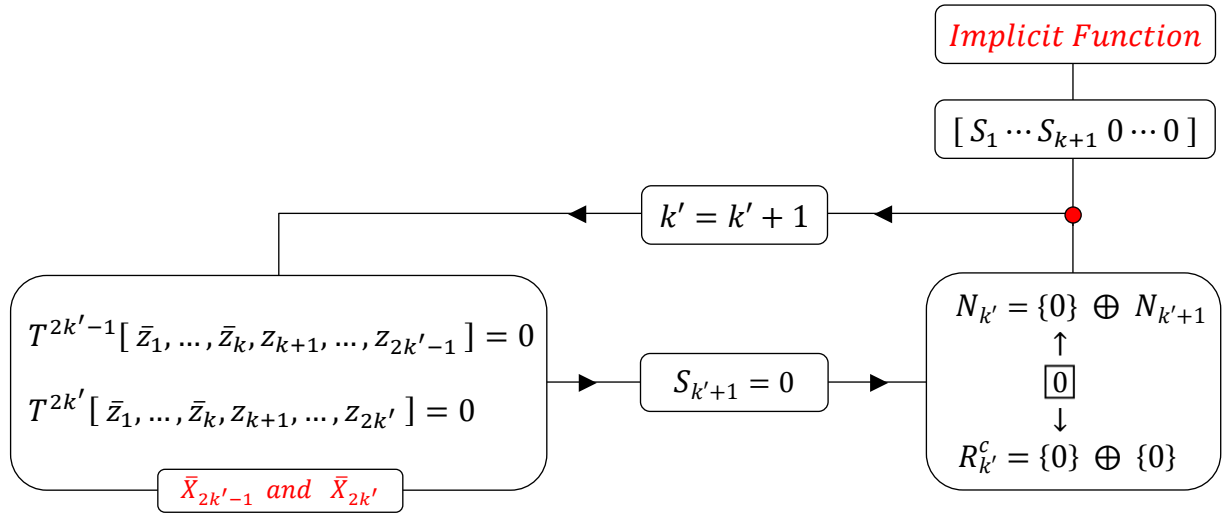

Figure 9 : The simplified iteration process with $[S_1 \cdots S_{k+1}]$ surjective.

Note that in figure 9, the previous iteration index $k$ as well as the leading coefficients $[\bar{z}_1, \ldots, \bar{z}_k]$ are fixed, whereas $k' = k+1, \ k+2, \ldots$ is used as a new running index.

Now, the main difference between figures 1 and 9 results from $S_{k'+1} = 0$ yielding trivial decomposition of $N_{k'}$ and $R_{k'}^c$ as well as trivial extension of the sum operator by $[S_1 \cdots S_{k+1}\ 0 \cdots 0]$. Hence, the filtration of $B$ stabilizes according to

$$B = N_0 \supset N_1 \supset \cdots \supset N_{k+1} = N_{k+2} = \cdots = N_{k'+1} \tag{55}$$

implying decomposition of $B$ and $\bar{B}$ simply given by

$$
\begin{array}{c}
\phantom{B = N_1^c \oplus N_2^c \oplus \cdots \oplus N_{k+1}^c \oplus} {=N_{k+2}^c} \phantom{\oplus \cdots \oplus} {=N_{k'+1}^c} \\
B = N_1^c \oplus N_2^c \oplus \cdots \oplus N_{k+1}^c \oplus \widetilde{\{0\}} \oplus \cdots \oplus \widetilde{\{0\}} \oplus N_{k'+1} \\
\uparrow \quad \uparrow \quad \quad \uparrow \quad \uparrow \quad \quad \uparrow \\
\boxed{S_1} \ \boxed{S_2} \quad \boxed{S_{k+1}} \ \boxed{0} \quad \boxed{0} \\
\downarrow \quad \downarrow \quad \quad \downarrow \quad \downarrow \quad \quad \downarrow \\
\bar{B} = R_1 \oplus R_2 \oplus \cdots \oplus R_{k+1} \oplus \{0\} \oplus \cdots \oplus \{0\} \\
\phantom{\bar{B} = R_1 \oplus R_2 \oplus \cdots \oplus R_{k+1} \oplus\ } {=R_{k+2}} \phantom{\oplus \cdots \oplus} {=R_{k'+1}}
\end{array}
\tag{56}
$$

Exploiting (55) and (56), it is not too difficult to see that the following properties hold concerning the correlation of derivatives between approximate/exact solution curves as well as the structure of $\bar{X}_\infty$.

## Corollary 2 - Derivatives and Structure of $\bar{X}_\infty$

Assume $k$-regularity of $[\bar{z}_1, \ldots, \bar{z}_k]$ and $k \geq 1$ to be minimal with respect to surjectivity of the sum operator $[S_1 \cdots S_{k+1}]$.

(i) For $l \geq 1$ we obtain

$$\pi_{k+l}(\bar{X}_{2k+l-1}) \supsetneq \pi_{k+l}(\bar{X}_{2k+l}) = \pi_{k+l}(\bar{X}_\infty),$$

i.e. $\pi_{k+l}(\bar{X}_i)$ stabilizes exactly at $i = 2k + l$, then agreeing with $\pi_{k+l}(\bar{X}_\infty)$.

(ii) The set $\bar{X}_{2k+l}, l \geq 1$ is explicitly parametrized by the subspace $N_1 \times \cdots \times N_k \times (N_{k+1})^l \subset B^{k+l}$ according to

$$[z_1, \ldots, z_k] = [\bar{z}_1, \ldots, \bar{z}_k]$$

---

$$z_{k+1} = P_1 + q_1 \qquad\qquad q_1 \in N_{k+1}, \ P_1 \in B$$

$$z_{k+2} = P_2(q_1) + q_2 \qquad\qquad q_2 \in N_{k+1}$$

$$\ldots \qquad\qquad \ldots$$

$$z_{k+l} = P_l(q_1, \ldots, q_{l-1}) + q_l \qquad\qquad q_l \in N_{k+1}$$

---

$$z_{k+l+1} = P_{l+1}(q_1, \ldots, q_l) + n_k \qquad\qquad n_k \in N_k$$

$$z_{k+l+2} = P_{l+2}(q_1, \ldots, q_l, n_k) + n_{k-1} \qquad\qquad n_{k-1} \in N_{k-1}$$

$$\ldots \qquad\qquad \ldots$$

$$z_{k+l+k} = P_{l+k}(q_1, \ldots, q_l, n_k, \ldots, n_2) + n_1 \qquad\qquad n_1 \in N_1$$

with operators $P_i(\cdot)$, $2 \leq i \leq l+k$ defined by the composition of multilinear mappings that are effectively constructed by the iteration from figure 1. If $B$ and $\bar{B}$ are of finite dimensions, e.g. $G: \mathbb{C}^n \to \mathbb{C}^m$ with $n > m$, then $P_i(\cdot)$ can be represented by polynomials of increasing degree with respect to chosen coordinates.

**Remarks 1)** Moving $l$ towards infinity allows by $\pi_{k+l}(\bar{X}_{2k+l}) = \pi_{k+l}(\bar{X}_\infty)$ and (ii) to calculate $\bar{X}_\infty$ in a constructive way, i.e. the subset of arc space $X_\infty$ with leading coefficients $[\bar{z}_1, \ldots, \bar{z}_k]$ reads

$$\bar{X}_\infty = \{\, (z_i)_{i \in \mathbb{N}} \mid [z_1, \ldots, z_k] = [\bar{z}_1, \ldots, \bar{z}_k], \ z_{k+l} = P_l(q_1, \ldots, q_{l-1}) + q_l, \ q_l \in N_{k+1}, \ l \in \mathbb{N} \,\}.$$

The inclusion $\bar{X}_\infty \subset \{\cdots\}$ follows straightforward by contradiction.

**2)** Again, due to $\pi_{k+l}(\bar{X}_{2k+l}) = \pi_{k+l}(\bar{X}_\infty)$, every element in $\bar{X}_{2k+l}$ ensures an exact solution curve with Taylor coefficients in $\bar{X}_\infty$ agreeing up to order $k+l$ with the given element from $\bar{X}_{2k+l}$. On contrary by $\pi_{k+l}(\bar{X}_{2k+l-1}) \supsetneq \pi_{k+l}(\bar{X}_{2k+l})$, the set $\pi_{k+l}(\bar{X}_{2k+l-1})$ contains elements that cannot be continued to $\bar{X}_\infty$.

In this sense, we obtain some sort of, so to say, weak Greenberg function $\bar{\beta}(\cdot)$ with respect to the subset $\bar{X}_\infty$ of arc space $X_\infty$ given by $\bar{\beta}(k+l) = 2k+l$, $l \geq 1$. In case of $l = 0$, the trivial identities $\pi_k(\bar{X}_{2k}) = \{\bar{z}_1, \ldots, \bar{z}_k\} = \pi_k(\bar{X}_\infty)$ hold true, suggesting $\bar{\beta}(k) := 2k$ and

$$\bar{\beta}(i) = k + i, \ \ i \geq k \quad \text{and} \quad \pi_i(\bar{X}_{\bar{\beta}(i)}) = \pi_i(\bar{X}_\infty).$$

With respect to arc space $X_\infty$, the function $\bar{\beta}(\cdot)$ defines, within the general setting of Banach spaces, a lower bound for a possibly existing true Greenberg function $\beta(i)$.

In addition, note that reparametrization of a solution curve $z(\varepsilon) = \varepsilon \bar{z}_1 + \cdots + 1/k! \, \varepsilon^k \bar{z}_k + \cdots$ by $\varepsilon \to \varepsilon^2$ delivers

$$z^*(\varepsilon) := z(\varepsilon^2) = \frac{1}{2}\varepsilon^2 \cdot \underbrace{(2\bar{z}_1)}_{=:z_2^*} + \cdots + \frac{1}{(2k)!}\varepsilon^{2k} \cdot \underbrace{\left(\frac{(2k)!}{k!}\bar{z}_k\right)}_{=:z_{2k}^*} + \cdots$$

showing, after some calculations, $(2k)$-regularity of new leading coefficients $[0, z_2^*, 0, \ldots, 0, z_{2k}^*]$ with corresponding weak Greenberg function $\bar{\beta}(i) = 2k + i$, $i \geq 2k$. Repeating this process by $\varepsilon \to \varepsilon^3, \varepsilon^4, \cdots$, we obtain a sequence of weak Greenberg functions $\bar{\beta}(i) = 3k + i$, $4k + i, \cdots$ derived from $k$-regularity of $[\bar{z}_1, \ldots, \bar{z}_k]$. Then, simply taking the maximum of these functions, a lower bound of a possibly existing true Greenberg function $\beta(i)$ of $X_\infty$ is constructed, as depicted in figure 10 in case of $k = 1, 2$ and 3.

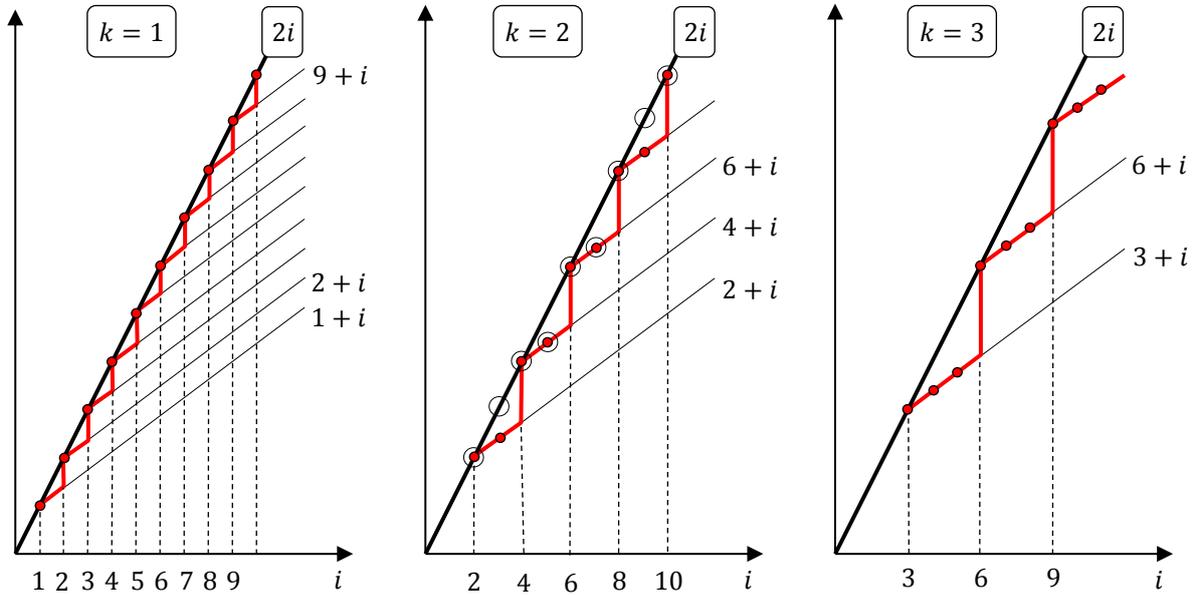

Figure 10 : Lower bound of Greenberg function $\beta(\cdot)$ derived from $[\bar{z}_1, \ldots, \bar{z}_k]$, $k = 1,2,3$.

In each case, the lower bound is defined by the red dotted points positioned on or below the straight line $2i$. For general $k \geq 1$, the leftmost point is given by $(k, 2k)$ on the line $2i$ followed by $k - 1$ points below the line $2i$ and a jump back to $2i$. In case of $k = 1$, no point below $2i$ occurs.

If the equation $G[z] = 0$, $G: B \to \bar{B}$ has several solution curves with regularity of degree $1 < k_1 < k_2 < \cdots$, then a superposition of the red marked lines in figure 10 occurs and for each $i \geq k_1$ the lower bound is defined by the maximum of corresponding red dotted points, finally yielding some sort of stepwise lower bound of $\beta(i)$. In the middle diagram of figure 10, the stepwise behavior is indicated by small circles in case of two solution curves with $k_1 = 2$ and $k_2 = 3$.

If we restrict to polynomials $G: \mathbb{C}^2 \to \mathbb{C}$, then the Greenberg function can explicitly be calculated according to [14], [23]. Other examples can be found in [21]. For instance, the Greenberg function of the monomial $G[x, y] = xy$ reads $\beta(i) = 2i$ showing that a lower bound cannot exceed the line $2i$ without imposing further assumptions on $G$.

Finally, in case of $G_0^1$ surjective ($k = 0$), the complete set of solutions is directly obtained by implicit function theorem and the Greenberg function $\beta(i)$ satisfies $\beta(i) = \bar{\beta}(i) = k + i = i$, i.e. the lower bound $\bar{\beta}(\cdot)$ already agrees with $\beta(i)$.

**3)** According to Remark 7) of Theorem 1, the iteration works quite well, when restricting to internal $x$-parametrization in case of $G: B = \mathbb{K} \times Y \to \bar{B}$. Hence, Corollary 2 can now be reinterpreted with respect to leading coefficients of the form

$$[\bar{z}_1, \cdots, \bar{z}_k] = [(1, \bar{y}_1)^T, (0, \bar{y}_2)^T, \ldots, (0, \bar{y}_k)^T]$$

with regularity of degree $k \geq 1$.

First, the iteration delivers all solutions within $\bar{X}_{2k+l}, l \geq 1$ that may still comprise coefficients $\bar{z}_i = (\bar{x}_i, \bar{y}_i)^T, i > k$ satisfying $\bar{x}_i \neq 0$, i.e. $\bar{X}_{2k+l}^0 \subset \bar{X}_{2k+l}$ and we have to restrict $\bar{X}_{2k+l}$ to elements with $\bar{x}_i = 0$. However, this can easily be established, yielding under consideration of decomposition (53) the following result with respect to $\pi_{k+l}(\bar{X}_{2k+l}^0)$.

The set $\pi_{k+l}(\bar{X}_{2k+l}^0), l \geq 1$ is parametrized by the subspace $(Y_{k+1})^l \subset Y^l$ such that

$$x_1 = 1, \ x_i = 0, \ i \geq 2$$
$$[y_1, \ldots, y_k] = [\bar{y}_1, \ldots, \bar{y}_k]$$
$$\text{---------------------------------------------------------}$$
$$y_{k+1} = p_1 + q_1 \qquad\qquad q_1 \in Y_{k+1}, \ p_1 \in Y$$
$$y_{k+2} = p_2(q_1) + q_2 \qquad\qquad q_2 \in Y_{k+1}$$
$$\ldots \qquad\qquad\qquad\qquad\qquad \ldots$$
$$y_{k+l} = p_l(q_1, \ldots, q_{l-1}) + q_l \quad q_l \in Y_{k+1}$$

with operators $p_i(\cdot), 2 \leq i \leq k + l$ defined again by the composition of multilinear mappings. The corresponding subset $\bar{X}_\infty^0 \subset X_\infty$ is given by

$$\bar{X}_\infty^0 = \left\{ \begin{pmatrix} x_i \\ y_i \end{pmatrix}_{i \in \mathbb{N}} \ \Big| \ x_1 = 1, x_{i \geq 2} = 0, [y_1, \ldots, y_k] = [\bar{y}_1, \ldots, \bar{y}_k], y_{k+l} = p_l(\cdot) + q_l, q_{l \geq 1} \in Y_{k+1} \right\}.$$

If $Y_{k+1} = \{0\}$, i.e. $N_{k+1} = \{(1, \bar{y}_1)^T\}$, then $\bar{X}_\infty^0$ contains exactly one element defined by $q_l = 0$, $l \geq 1$. The corresponding power series represents an exact analytic solution curve given by internal $x$-parametrization. Note also that the finiteness of $\bar{X}_\infty^0$ is in sharp contrast to $\bar{X}_\infty$, where every element gives rise to an infinity of other elements within $\bar{X}_\infty$ obtained by arbitrary $\varepsilon$-reparametrization. The single element in $\bar{X}_\infty^0$ is simply sorted out of $\bar{X}_\infty$ by posing the condition of being internally parametrized by $x$.

The case $Y_{k+1} = \{0\}$ arises when considering $G: \mathbb{K} \times \mathbb{K}^n \to \mathbb{K}^n$. Then, there may exist several isolated solution curves through the singularity at 0, parametrized internally by $x$ with leading coefficients of the form $\left[(1, \bar{y}_1^1)^T, \ldots, (0, \bar{y}_{k_1}^1)^T\right], \left[(1, \bar{y}_1^2)^T, \ldots, (0, \bar{y}_{k_2}^2)^T\right], \ldots$, showing different degrees of regularity $k_1, k_2, \ldots$, as well as different sets $\bar{X}_\infty^{0,1}, \bar{X}_\infty^{0,2}, \ldots$ (each composed of one element).

If $Y_{k+1} \neq \{0\}$ and $dim(Y_{k+1})$ is finite, then the solution curves are embedded within surfaces of solutions with dimension $1 + dim(Y_{k+1}) \geq 2$. In this case, $\bar{X}_\infty^0$ also contains an infinity of elements caused by different curves on the surface.

Concerning the weak Greenberg function with respect to the subset $\bar{X}_\infty^0$ of arc space $X_\infty$, we obtain again $\bar{\beta}(i) = k + i, \ i \geq k$.

Up to now, we mainly concentrated on pure existence of solution curves of $G[z] = 0$, thus describing the underlying solutions of the system of undetermined coefficients $T^1 = 0, T^2 = 0, \ldots$ as precisely as possible.

The main result is given by the family of solution curves from Theorem 1

$$z(\varepsilon,\bar{n}) = \varepsilon \bar{z}_1 + \cdots + \frac{1}{k!}\varepsilon^k \bar{z}_k + \frac{1}{(k+1)!}\varepsilon^{k+1} z_{k+1}(\varepsilon,\bar{n}) + \cdots + \frac{1}{(2k+1)!}\varepsilon^{2k+1} z_{2k+1}(\varepsilon,\bar{n}) \qquad (57)$$

parametrized with respect to $(\varepsilon,\bar{n}) = (\varepsilon, n_1, \ldots, n_k, n_{k+1}) \in \mathbb{R} \times N_1 \times \cdots \times N_{k+1}$. Obviously, the family shows a strong redundancy due to possible reparametrization (e.g. $\varepsilon \to \varepsilon + \varepsilon^{k+1}$), yielding same solution orbits of $G[z] = 0$ but different elements within $\bar{X}_{2k}$ and $\bar{X}_{2k+1}$. In the next step, we show how to reduce the parametrization within the family $z(\varepsilon,\bar{n})$ considerably without losing solutions of $G[z] = 0$.

First, note that there exists a minimal number $l$ with $\bar{z}_l \neq 0$, $1 \leq l \leq k$, because of $[S_1 \cdots S_{k+1}] = [0 \cdots 0]$ in case of $\bar{z}_1 = \cdots \bar{z}_k = 0$, which contradicts our assumption of $[S_1 \cdots S_{k+1}]$ to be surjective. In addition, it is not difficult to see that $\bar{z}_l$ is contained in every null space of the filtration, i.e. also in the smallest null space $N_{k+1}$, allowing us to refine decomposition (43) according to

$$B = N_1^c \oplus \cdots \oplus N_{k+1}^c \oplus \underbrace{\overbrace{N_{k+1}/\{\bar{z}_l\} \oplus \{\bar{z}_l\}}^{=N_{k+1}}}_{=:P} \qquad (58)$$

Here $\{\bar{z}_l\}$ denotes the space in $B$ spanned by $\bar{z}_l$ and $P$ labels a direct complement of $\{\bar{z}_l\}$ within $N_{k+1}$. Now, by restricting the family of solution curves $z(\varepsilon, n_1, \ldots, n_k, n_{k+1})$ to the subspace $P$ by

$$z_{red}(\varepsilon, p) := z(\varepsilon, 0, \ldots, 0, p), \quad p \in P \qquad (59)$$

we will see that essentially no solutions were lost.

Under consideration of (57) and Theorem 1 (iii), the map $z_{red}(\varepsilon, p)$ can be written in detail by

$$
\begin{aligned}
z_{red}(\varepsilon, p) &= \tfrac{1}{l!}\varepsilon^l \bar{z}_l + \cdots + \tfrac{1}{k!}\varepsilon^k \bar{z}_k + \tfrac{1}{(k+1)!}\varepsilon^{k+1} z_{k+1}(\varepsilon, p) + \cdots + \tfrac{1}{(2k+1)!}\varepsilon^{2k+1} z_{2k+1}(\varepsilon, p) \\
&= \underbrace{\tfrac{1}{l!}\varepsilon^l \bar{z}_l + \cdots + \tfrac{1}{k!}\varepsilon^k \bar{z}_k}_{} \\
&\quad + \left[\tfrac{1}{(2k+1)!}\varepsilon^{2k+1} \cdots \tfrac{1}{(k+1)!}\varepsilon^{k+1}\right] \cdot \left\{\hat{I}_{k+1} + \boxed{\widehat{M}_{k+1}} \cdot \left[\begin{pmatrix} n_1^c \\ \vdots \\ n_{k+1}^c \end{pmatrix}(\varepsilon, p) + \boxed{\widehat{L}_{k+1}} \cdot \begin{pmatrix} 0 \\ \vdots \\ p \end{pmatrix}\right]\right\} \\[-2pt]
&\qquad \qquad \qquad \qquad \qquad \qquad \qquad \qquad \qquad =: z_0(\varepsilon)
\end{aligned}
\qquad (60)
$$

$$
\begin{aligned}
&= \tfrac{1}{l!}\varepsilon^l \bar{z}_l + \cdots + \tfrac{1}{k!}\varepsilon^k \bar{z}_k + \left[\tfrac{1}{(2k+1)!}\varepsilon^{2k+1} \cdots \tfrac{1}{(k+1)!}\varepsilon^{k+1}\right] \cdot \hat{I}_{k+1} \\
&\quad + \underbrace{\left[\tfrac{1}{(2k+1)!}\varepsilon^{2k+1} \cdots \tfrac{1}{(k+1)!}\varepsilon^{k+1}\right] \cdot \left\{\boxed{\widehat{M}_{k+1}} \cdot \left[\begin{pmatrix} n_1^c \\ \vdots \\ n_{k+1}^c \end{pmatrix}(\varepsilon, p) + \boxed{\widehat{L}_{k+1}} \cdot \begin{pmatrix} 0 \\ \vdots \\ p \end{pmatrix}\right]\right\}}_{=: z_1(\varepsilon, p)}
\end{aligned}
$$

$$=: z_0(\varepsilon) + z_1(\varepsilon, p).$$

Next, taking account of Theorem 1 (ii), the coefficients $[\bar{z}_1, \ldots, \bar{z}_k, \hat{I}_{k+1,1}, \ldots, \hat{I}_{k+1,k+1}]$ of $z_0(\varepsilon)$ define an element in $\bar{X}_{2k+1}$ which determines some sort of orientation line of the solution curves $z_{red}(\varepsilon, p) = z_0(\varepsilon) + z_1(\varepsilon, p)$, as qualitatively indicated in figure 11 below.

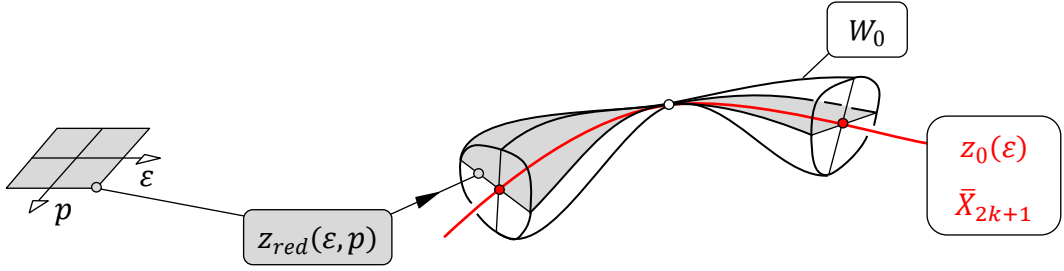

Figure 11 : The smooth and regular solution manifold $z_{red}(\varepsilon, p)$ in $W_0$.

Then choose $p \in P$ with $|p| < c, c > 0$ arbitrary large but fixed and $\varepsilon$ sufficiently small.

**Corollary 3 - Regularity and Uniqueness**

(i) The map $z_{red}(\varepsilon, p)$ defines a regular Banach manifold in $B$ for $\varepsilon \neq 0$.

(ii) There exists an open neighborhood $W_0$ of $z_0(\varepsilon)$, $\varepsilon \neq 0$ with the property that all solutions of $G[z] = 0$ in $W_0$ are given by $z_{red}(\varepsilon, p)$.

(iii) The solutions $z(\varepsilon, n_1, \dots, n_k, \alpha \bar{z}_l + p)$, $\alpha \in \mathbb{K}$, $p \in P$ from Theorem 1 (i) are contained in the manifold $z_{red}(\varepsilon, p)$ for $n_1, \dots, n_k$ and $\alpha$ sufficiently small (uniformly in $\varepsilon$).

If the dimension $dim(P)$ of the subspace $P$ is finite, then the dimension of the manifold equals $1 + dim(P)$. Figure 11 shows a constellation with $dim(P) = 1$, whereas $dim(P) = 0$ is indicated within figure 12 below, yielding exactly one solution orbit within $W_0$.

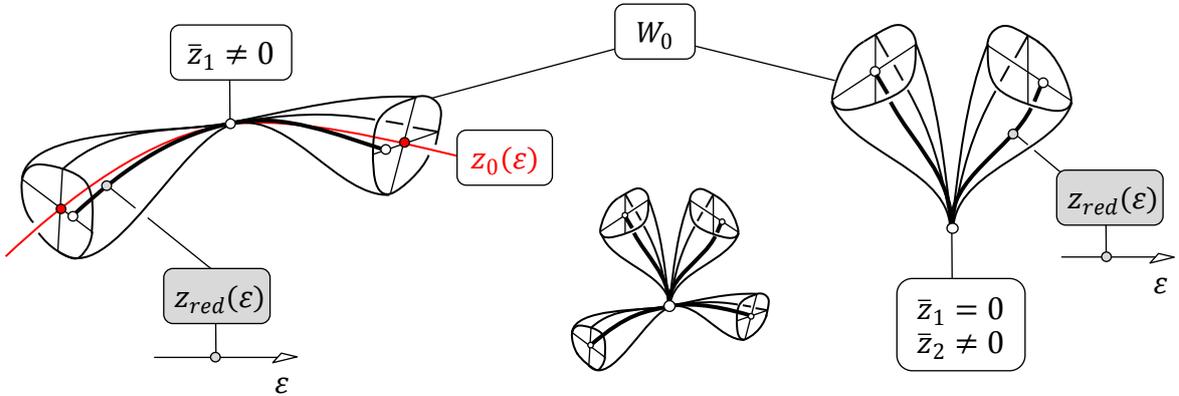

Figure 12 : The case $dim(P) = 0$ with $\bar{z}_1 \neq 0$ (left) and $\bar{z}_1 = 0, \bar{z}_2 \neq 0$ (right).

Note that in general the red center line $z_0(\varepsilon)$ with coefficients from $\bar{X}_{2k+1}$ does not agree with the solution orbit as indicated in the left diagram. In the middle diagram, the typical constellation of pointed wedges of uniqueness near a bifurcation point of two solution curves is shown.

If $G$ is given by a map of the form $G: \mathbb{R}^{n+q} \to \mathbb{R}^n$, then $q = 1 + dim(P)$ due to bijectivity of the sum operator $[S_1 \cdots S_{k+1}]$ with respect to $N_1^c \oplus \cdots \oplus N_{k+1}^c$ and $\bar{B} = \mathbb{R}^n$. Hence, figure 11 is obtained in case of $G: \mathbb{R}^{n+2} \to \mathbb{R}^n$, whereas figure 12 occurs in case of $G: \mathbb{R}^{n+1} \to \mathbb{R}^n$.

We sketch the proofs of (i)-(iii) within the following remarks.

**Remarks 1)** The regularity of the manifold, i.e. $N[z'_{red}(\varepsilon, p)] = \{0\}$, $\varepsilon \neq 0$, results from (58) and (60) along the following lines. Consider the equation

$$z'_{red}(\varepsilon, p) \cdot \begin{pmatrix} \bar{\varepsilon} \\ \bar{p} \end{pmatrix} = 0$$

$$\Leftrightarrow \left[ \tfrac{1}{(l-1)!} \varepsilon^{l-1} \cdot \bar{z}_l + O(\varepsilon^l) \right] \cdot \bar{\varepsilon}$$

$$+ \left[ \tfrac{1}{(2k+1)!} \varepsilon^{2k+1} \cdots \tfrac{1}{(k+1)!} \varepsilon^{k+1} \right] \cdot \boxed{\widehat{M}_{k+1}} \cdot \left[ \begin{pmatrix} n^c_{1,p} \\ \vdots \\ n^c_{k+1,p} \end{pmatrix}(\varepsilon, p) + \boxed{\widehat{L}_{k+1}} \cdot \begin{pmatrix} 0 \\ \vdots \\ 1 \end{pmatrix} \right] \cdot \bar{p} = 0$$

Here $n^c_{i,p}(\varepsilon, p) \in L[P, N^c_i]$ denotes the Fréchet derivative of $n^c_i(\varepsilon, p)$ with respect to $p \in P$. Then, with $\varepsilon \neq 0$ sufficiently small, we obtain $\bar{\varepsilon} = 0$ due to $\bar{z}_l \neq 0$, and it suffices to discuss the remaining equation which turns out to be of the form

$$\tfrac{1}{(k+1)!} \varepsilon^{k+1} \cdot \left[ n^c_{k+1,p}(\varepsilon, p) \cdot \bar{p} + \bar{p} \right] + O(|\varepsilon|^{k+2}) = 0$$

by using some properties of $\widehat{M}_{k+1}$ and $\widehat{L}_{k+1} \in GL[B^{k+1}, B^{k+1}]$. In addition, $N^c_{k+1} \cap P = \{0\}$ by (58), yielding $\bar{p} = 0$ and the regularity of the manifold in (i) is shown with respect to $\varepsilon \neq 0$.

**2)** The wedge-like neighborhood $W_0 \subset B$ of $z_0(\varepsilon)$ from (ii) is defined using the image of the map

$$A: \underbrace{(\varepsilon, n^c_1, \ldots, n^c_{k+1}, p)}_{\in \mathbb{K} \times N^c_1 \times \cdots \times N^c_{k+1} \times P} \rightarrow z_0(\varepsilon) + \left[ \tfrac{1}{(2k+1)!} \varepsilon^{2k+1} \cdots \tfrac{1}{(k+1)!} \varepsilon^{k+1} \right] \cdot \boxed{\widehat{M}_{k+1}} \cdot \left[ \begin{pmatrix} n^c_1 \\ \vdots \\ n^c_{k+1} \end{pmatrix} + \boxed{\widehat{L}_{k+1}} \cdot \begin{pmatrix} 0 \\ \vdots \\ p \end{pmatrix} \right] \tag{61}$$

The first summand $z_0(\varepsilon)$ defines the center line of the map with leading term given by $\bar{z}_l$. The second summand defines a linear map with respect to the variables $[n^c_1, \ldots, n^c_{k+1}, p]$ and range given by $N^c_1 \oplus \cdots \oplus N^c_{k+1} \oplus P$ which defines by itself a direct complement to $\bar{z}_l$ within $B$, i.e.

$$B = \{\bar{z}_l\} \oplus \{ N^c_{k+1} \oplus \cdots \oplus N^c_{k+1} \oplus P \}$$

by (58). Hence, at every point of the center line $z_0(\varepsilon)$ a direct complement of $\bar{z}_l$ within $B$ is attached.

Then, by restriction of $[n^c_1, \ldots, n^c_{k+1}, p]$ to a bounded region comprising zero, the image of the second summand is shrinking to zero as $\varepsilon \rightarrow 0$, thereby creating the open and wedge-like neighborhood $W_0$ of $z_0(\varepsilon)$ for $\varepsilon \neq 0$. In some more detail, the elements from $N^c_{k+1}$ and $P$ are multiplied by $\varepsilon^{k+1}$ in (61), implying decrease of $W_0$ in the directions of $N^c_{k+1}$ and $P$ by order of $O(|\varepsilon|^{k+1})$. The remaining subspaces $N^c_l, 1 \leq l \leq k$ show faster shrinking at least by order of $O(|\varepsilon|^{k+2})$.

The uniqueness of the solutions $z_{red}(\varepsilon, p)$ in $W_0$ can now be seen as follows. Apply the implicit function theorem to the remainder equation

$$T^{2k+1}\left[ \widehat{I}_{k+1} + \boxed{\widehat{M}_{k+1}} \cdot \left[ \begin{pmatrix} n^c_1 \\ \vdots \\ n^c_{k+1} \end{pmatrix} + \boxed{\widehat{L}_{k+1}} \cdot \begin{pmatrix} 0 \\ \vdots \\ p \end{pmatrix} \right] \right] + O(|\varepsilon|) = 0 \tag{62}$$

choosing $[n^c_1, \ldots, n^c_{k+1}, p]$ from the bounded region comprising zero from above. Locally, we obtain uniqueness, carrying directly over to $W_0 \subset B$ and ascertaining (ii). Finally, show by continuity that the complete family of solution curves $z(\varepsilon, n_1, \ldots, n_k, \alpha\bar{z}_l + p), \alpha \in \mathbb{K}, p \in P$ from Theorem 1 (i) is entering $W_0$ independent of $\varepsilon$ for $n_1, \ldots, n_k, \alpha$ sufficiently small, thus implying (iii).

## 4. An application to Newton-Polygons and the Milnor Number

In this section consider $G : \mathbb{K}^2 \to \mathbb{K}$, $\mathbb{K} = \mathbb{R}, \mathbb{C}$ analytic near $z = (x,y) = (0,0) \in \mathbb{K}^2$ with

$$G[x,y] = \sum_{\alpha,\beta=0}^{\infty} c_{\alpha\beta} \cdot x^\alpha \cdot y^\beta \quad \text{and} \quad G[0,0] = 0. \qquad (63)$$

Possible solutions of the form $x \equiv 0$ or $y \equiv 0$ are supposed to be split off, implying a convenient Newton-polygon as indicated in figure 13 left.

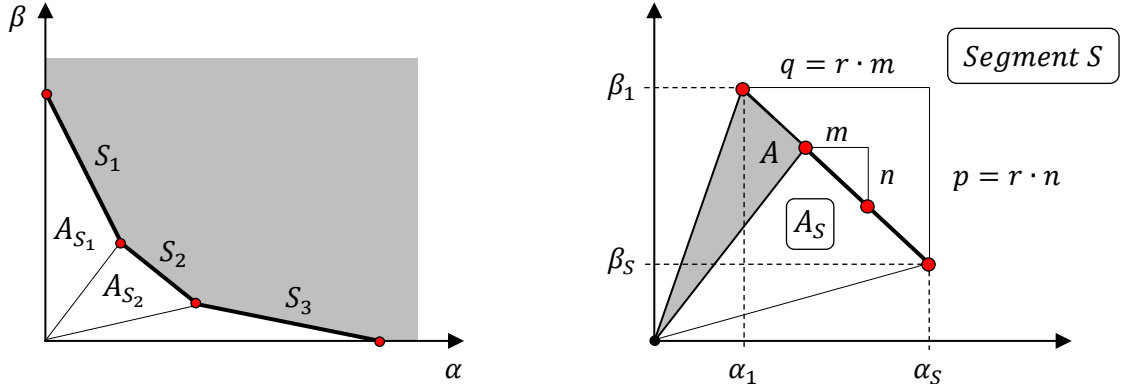

Figure 13 : Convenient Newton-polygon (left) with one segment $S$ in detail (right).

The total number of segments $S_i$ will be denoted by $\tau \geq 1$. Every segment is individually characterized by integer data

$$p = r \cdot n, \quad q = r \cdot m, \quad r \geq 1, \quad gcd(n,m) = 1 \quad \text{and} \quad N = n \cdot \alpha_S + m \cdot \beta_S \qquad (64)$$

as shown in figure 13 right. $N$ denotes the degree of the quasi-homogenous polynomial associated with the given segment. Area $A$ is defined by the grey area in figure 13 right, and $A_S$ represents the complete area defined by the triangle of origin and the two end points of the segment. Then blowing up by

$$x = \varepsilon^n \cdot v, \quad y = \varepsilon^m \cdot w$$

equations (63) and (64) imply

$$G[\varepsilon^n v, \varepsilon^m w] = \varepsilon^N \cdot \{v^{\alpha_1} \cdot w^{\beta_S} \cdot [c_{\alpha_1\beta_1}(w^n)^r + c_{\alpha_2\beta_2}(w^n)^{r-1}(v^m)^1 + \cdots + c_{\alpha_S\beta_S}(v^m)^r] + O(|\varepsilon|)\}$$

$$= \varepsilon^N \cdot \{v^{\alpha_1} \cdot w^{\beta_S} \cdot \phi(v^m, w^n) + O(|\varepsilon|)\} \qquad (65)$$

with $\phi(.,.)$ homogenous of degree $r \geq 1$. Now, if $\mathbb{K} = \mathbb{C}$, complete factorization holds with respect to the homogenous polynomial $\phi(.,.)$, yielding the equation

$$H[\varepsilon,v,w] := v^{\alpha_1} \cdot w^{\beta_S} \cdot (a_1 v^m + b_1 w^n) \cdot \cdots \cdot (a_r v^m + b_r w^n) + O(|\varepsilon|) = 0 \qquad (66)$$

with $a_i, b_i \neq 0$ and identical factors not excluded. If $\mathbb{K} = \mathbb{R}$, we end up with partial factorization of $\phi(.,.)$ by

$$H[\varepsilon,v,w] := v^{\alpha_1} \cdot w^{\beta_S} \cdot (a_1 v^m + b_1 w^n) \cdot \cdots \cdot (a_t v^m + b_t w^n) \cdot R(v^m, w^n) + O(|\varepsilon|) = 0 \quad (67)$$

with $0 \leq t \leq r$ and remainder polynomial $R(.,.)$. Next, choose the maximum of $m$ and $n$, say $n$ for definiteness and assume $a_1 v^m + b_1 w^n$ to be a simple factor with multiplicity 1 within (66) or (67). Further, let $v_0 \neq 0$, $w_0 \neq 0$ be an arbitrary solution of $a_1 v^m + b_1 w^n = 0$, yielding $H[0, v_0, w_0] = 0$ and regularity by $H_v[0, v_0, w_0] \neq 0$, finally implying the existence of a smooth solution curve $v(\varepsilon)$ with $v(0) = v_0$ and $H[\varepsilon, v(\varepsilon), w_0] = 0$.

In addition, by uniqueness and reparametrization, it is easily seen that choosing other base solutions $v_0 \neq 0$, $w_0 \neq 0$ will result in the same solution orbit. Hence, every factor within (66) or (67) of multiplicity 1 gives rise to unique solutions of $G[x,y] = 0$ according to

$$\begin{pmatrix} x \\ y \end{pmatrix} = \varepsilon^m \cdot \begin{pmatrix} 0 \\ w_0 \end{pmatrix} + \varepsilon^n \cdot \begin{pmatrix} v_0 \\ 0 \end{pmatrix} + \begin{pmatrix} O(|\varepsilon|^{n+1}) \\ 0 \end{pmatrix} \tag{68}$$

$$=: \frac{1}{m!} \varepsilon^m \cdot \bar{z}_m + \frac{1}{n!} \varepsilon^n \cdot \bar{z}_n + O(|\varepsilon|^{n+1}).$$

Generically, in case of $\mathbb{K} = \mathbb{C}$, the multiplicity of each linear factor equals 1, yielding $r$ different solution curves of $G[z] = 0$ derived from the given segment $S$. In case of $\mathbb{K} = \mathbb{R}$, the number of different solution curves decreases correspondingly.

Up to now, the calculations are standard with respect to Newton-polygons [6]. Next, we aim to calculate the $k$-degree of regularity of the solution curves (68). For this purpose, the expansion (68) is plugged into the iteration from figure 1, looking for the first value of $k$ with $S_{k+1} \neq [0\ 0]$, hence implying surjectivity of the sum operator. Then, under consideration of the quasi-homogenous degree $N$ from (64) and the areas $A$ and $A_S$ from figure 13 right, we obtain by basic calculation the following results.

**Corollary 4 - Newton-polygons and Milnor number**

(i) Every linear factor of multiplicity 1 of segment $S$ defines a solution curve by (68) with $k$-regularity of minimal degree

$$k = N - \max(m, n) = 2A - \max(m, n). \tag{69}$$

(ii) If the homogenous polynomial $\phi(.,.)$ of segment $S$ factorizes completely with multiplicities 1, then $r$ different solution curves arise, each having minimal $k$-degree given by (69) and summing up to

$$r \cdot k = 2A_S - \max(p, q). \tag{70}$$

(iii) If the homogenous polynomials $\phi_i(.,.)$, $i = 1, \ldots, \tau$ of all segments (with data $m_i, n_i, r_i$) factorize completely with multiplicities 1, then $r_1 + \cdots + r_\tau$ different solution curves exist and the Milnor number $\mu$ of the singularity can be calculated by

$$\mu = \sum_{i=1}^{\tau} r_i \cdot k_i - ord(G) + 1. \tag{71}$$

Here $k_i = N_i - \max(m_i, n_i)$, $i = 1, \ldots, \tau$ by (69) and $ord(G)$ denotes lowest degree of $G$.

The relation between $k$-degrees of solution curves and the Milnor number $\mu$ of the singularity in (iii) follows directly from (i), (ii) and Kouchnirenko's planar theorem [1], [17].

**Remark:** If $G$ in (63) satisfies $c_{10} \neq 0$ or $c_{01} \neq 0$, i.e. if the Newton-polygon is composed of a single segment with endpoint at $\alpha = 1$ or $\beta = 1$, then $S_1 = G_0^1 = [c_{10}\ c_{01}] \neq [0\ 0]$ is already surjective and the implicit function theorem can be applied without entering the iteration in figure 1. In this sense, minimal $k$-degree of the solution curve is given by $k = 0$, agreeing with (69).

**Example:** Let us look at the simple constellation of a power series (63) with Newton-polygon given by a homogenous polynomial of degree $N$ as indicated in figure 14.

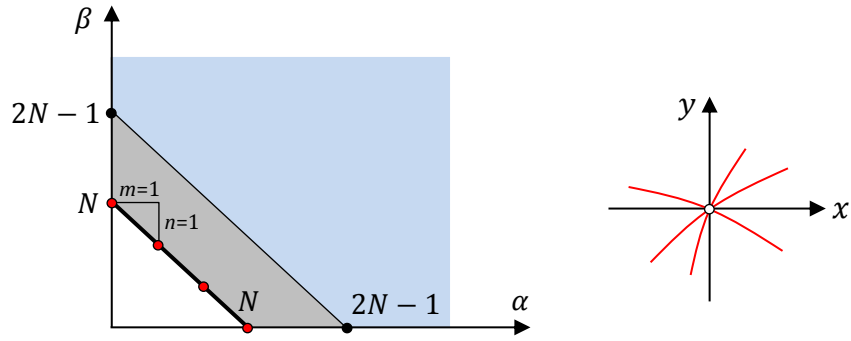

Figure 14 : Newton-polygon of homogenous degree $N$ (left) with typical solution curves (right).

Then, minimal $k$-degree of solution curves originating from linear factors with multiplicity 1 reads by (69)

$$k = N - 1$$

and if additionally $N$ different linear factors exist, the Milnor number $\mu$ is given by

$$\mu = N \cdot (N - 1) - N + 1 = (N - 1)^2$$

according to (71). Finally, by the stability result of Corollary 1 in section 3, monomials of order $\geq 2k + 1 = 2N - 1$ (within blue area) have no influence on the derivatives $[\bar{z}_1, \ldots, \bar{z}_{N-1}]$ of the solution curves. A typical bifurcation diagram in $\mathbb{R}^2$ is shown on the right hand side in figure 14. Obviously, in case of $N = 1$, we obtain $k = 0$ and $\mu = 0$.

In the next section, we aim to relax the assumption of simple zeros with multiplicity 1 within Corollary 4 by using a generalization of Kouchnirenko's planar theorem proved in [6].

### 4.1 Higher multiplicities and strict transforms

Let us name the segments of the Newton-polygon with $p > q$ and $p < q$ upper and lower segments respectively. If a segment with $p = q$ exists, it is called the neutral segment as shown in figure 15 left.

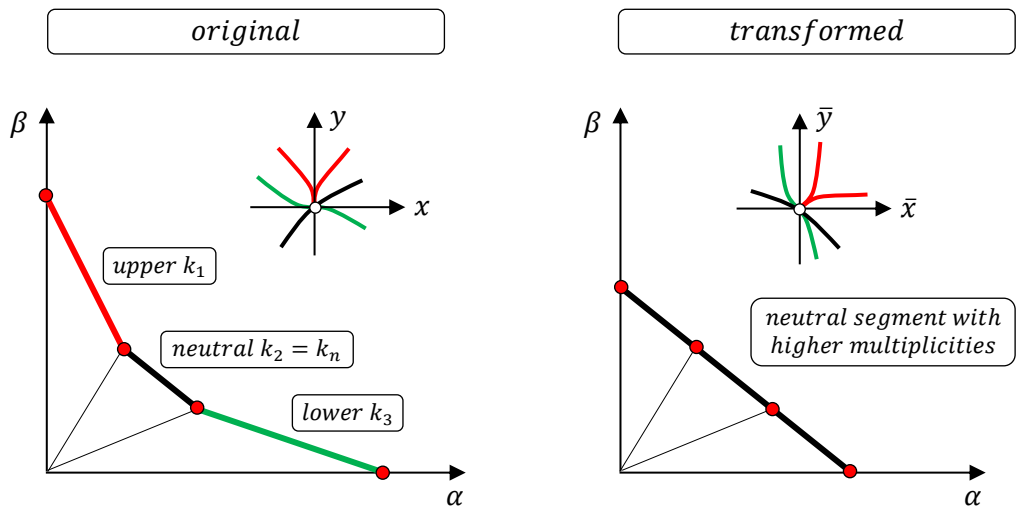

Figure 15 : Newton-polygon and solution curves; original (left) and transformed (right).

Then, from Corollary 4 (i), it is easily seen that $k$-degree of simple zeros achieves its minimal value at the neutral segment with increasing behavior when moving to upper and lower segments. More precisely,

$$k_1 > k_2 > \cdots > k_{n-1} \geq k_n \quad \text{and} \quad k_n \leq k_{n+1} < \cdots k_{\tau-1} < k_\tau \tag{72}$$

where equality only occurs, if the neutral segment has an endpoint with $\alpha = 1$ or $\beta = 1$.

Geometrically, in $\mathbb{R}^2$ upper segments imply solution curves of the form (68) with $m < n$ tangentially touching the $y$-axis, whereas solution curves belonging to lower segments are tangentially touching the $x$-axis, as indicated in the $(x, y)$-coordinate system on the left hand side within figure 15.

Now, when performing a diffeomorphic coordinate transformation of $(x, y)$ variables, e.g. a small rotation, then the solution curves generically leave the coordinate axes, as depicted in figure 15 top right, thereby turning into solution curves that belong to the neutral segment of the transformed Newton-polygon. In this sense, the neutral segment typically absorbs the solution curves of upper and lower segments under transformation.

On the other hand, high $k$-degree of solution curves belonging to upper and lower segments should not be affected by the transformation (cf. Remark 5 of Theorem 1) with the consequence that the transformed solution curves cannot appear as simple zeros within the neutral segment of the transformed system but only as zeros with higher multiplicities.

This is the first motivation that it should be possible to calculate $k$-degree of regularity from data of the Newton-polygon also in case of higher multiplicities. Secondly, the Milnor number $\mu$ is invariant with respect to diffeomorphic coordinate transformations and one may ask, whether formula (71) remains structurally valid if linear factors with higher multiplicities arise within the homogenous polynomials $\phi_i(.,.)$, $i = 1, \ldots, \tau$ of the segments?

Now, with respect to this investigation, assume $\mathbb{K} = \mathbb{C}$ for simplicity, ensuring factorization (66) of segment S according to

$$H[\,\varepsilon, v, w\,] := v^{\alpha_1} \cdot w^{\beta_s} \cdot (a_1 v^m + b_1 w^n)^{d_1} \cdot \cdots \cdot (a_\sigma v^m + b_\sigma w^n)^{d_\sigma} + O(\,|\varepsilon|\,) = 0$$

with $\sigma$ different factors of multiplicities $d_1, \ldots, d_\sigma$ satisfying $d_1 + \cdots + d_\sigma = r$. Then, move the base solution $v_0 \neq 0$, belonging to the first factor $(a_1 v^m + b_1 w^n)^{d_1}$ to the origin according to $v = v_0 + \bar{v}$ and fix the variable $w$ to its base solution by $w = w_0$, ending up with the transformed equation

$$\bar{H}[\,\varepsilon, \bar{v}\,] := H[\,\varepsilon, v_0 + \bar{v}, w_0\,] = c_1 \bar{v}^{d_1} + c_2 \varepsilon^{e_1} + O(\,|\bar{v}|^{d_1+1}\,) + O(\,|\varepsilon|^{e_1+1}\,) = 0$$

and $c_1 \neq 0$. Further assume $c_2 \neq 0$. The power series $\bar{H}[\varepsilon, \bar{v}]$ is typically called the strict transform belonging to the linear factor $(a_1 v^m + b_1 w^n)^{d_1}$ with corresponding Newton-polygon as depicted in figure 16 right.

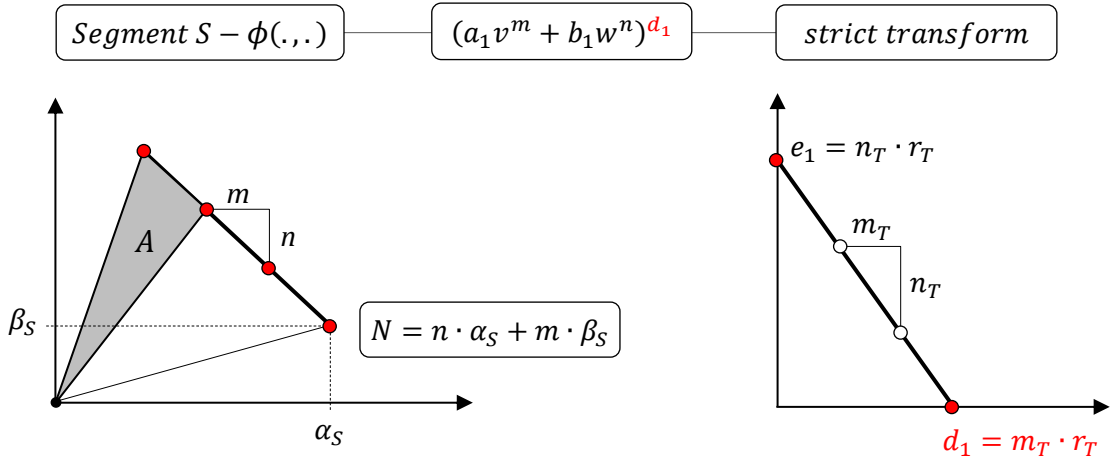

Figure 16 : Segment S with linear factor of multiplicity $d_1$ and corresponding strict transform.

The Newton polygon of the strict transform is characterized by one segment with end points on the axes (due to generic assumption $c_2 \neq 0$) and possible integer points in between, defining $(m_T, n_T)$ and $\gcd(d_1, e_1) = r_T$ as usual. Note that $d_1$ is restricted by $r$, whereas $e_1 \geq 1$ can assume arbitrary values.

Then, by direct calculation, exactly $r_T$ different solution branches result from the Newton polygon of the strict transform, implying $r_T$ different solution branches of the original equation $G[x, y] = 0$ derived from the factor $(a_1 v^m + b_1 w^n)^{d_1}$ of segment $S$.

Now, with respect to $k$-degree of regularity of these solution branches as well as the Milnor number $\mu$ of the singularity, we expect the following results to hold true.

- Minimal $k$-degree of regularity of each of the $r_T$ solution branches is given by

$$k = m_T \cdot [\, N - max(m, n) \,] + n_T \cdot [d_1 - 1]. \qquad (73)$$

- The Milnor number $\mu$ can be calculated by $k$-degree of all solution branches according to

$$\mu = \sum_{i=1}^{\tau} \sum_{j=1}^{\sigma_i} r_{ij} \cdot k_{ij} - ord(G) + 1. \qquad (74)$$

Here $r_{ij}$ denotes the number of different solution branches derived from the strict transform belonging to segment $i, 1 \leq i \leq \tau$ and linear factor $j, 1 \leq j \leq \sigma_i$, whereas $k_{ij}$ denotes the corresponding degree of regularity according to (73).

Considering (73), a zero of the Newton-polygon with arbitrary multiplicity gives rise to $r_{ij}$ solution branches, each having minimal $k$-degree by (73) that is calculated by geometric data $(N, m, n)$ of the segment, by the multiplicity $d_1$ of the zero, and finally by geometric data $(m_T, n_T)$ of the strict transform.

Formula (74) follows directly from (73) by use of Corollary 4 and the generalization of Kouchnirenko's planar theorem given in [6]. It remains to prove (73), possibly by inspection of the iteration from figure 1.

*Matthias Stiefenhofer*

*University of Applied Sciences*

*87435 Kempten (Germany)*

*matthias.stiefenhofer@hs-kempten.de*